\documentclass[12pt]{article}

\usepackage[utf8]{inputenc}
\usepackage[english]{babel}
\usepackage[tbtags]{amsmath}
\usepackage{physics}
\usepackage{amsfonts,amssymb,mathrsfs,amscd}
\usepackage{hyperref}
\usepackage{geometry}
\usepackage{titling}
\usepackage{amsthm}
\usepackage[T1,T2A]{fontenc}

\makeatletter
\gdef\No{{\select@language{english}\textnumero}}
\makeatother

\numberwithin{equation}{section}
\theoremstyle{plain}
\newtheorem{theorem}{Theorem}
\newtheorem{lemma}{Lemma}[section]
\newtheorem{propos}{Proposition}
\newtheorem{cons}{Corollary}
\theoremstyle{definition}
\newtheorem{definition}{Definition}
\newtheorem{remark}{Remark}[section]

\def\NN{\mathbb N}

\def\O{\mathcal{O}}

\def\Z{\mathbb Z}
\def\Q{\mathbb Q}
\def\V{\mathbb V}
\def\G{\mathbb G}
\def\Homgr{\underline{\mathrm{Hom}}^{\mathrm{gr}}}
\def\sgn{\mathrm{sgn}}
\def\Spec{\mathrm{Spec}}

\def\Nil{\mathrm{Nil}}
\def\Ladd{\mathcal{L}^{n}(\mathbb{G}_{a})}
\def\Lmult{\mathcal{L}^{n}(\mathbb{G}_{m})}
\begin{document}

\title{Polymultiplicative maps associated with the algebra of Iterated Laurent series and the higher-dimensional Contou-Carrere Symbol}
\author{Vladislav Levashev}
\date{}

\maketitle
\vspace*{-2cm}

\begin{abstract}
We study functorial polymultiplicative maps from the multiplicative group of the algebra of $n$-times iterated Laurent series over a commutative ring in $n+1$ variables into the multiplicative group of the ring. It is proven that if such a map is invariant under continuous automorphisms of this algebra, then it coincides, up to a sign, with an integer power of the $n$-dimensional Contou-Carrère symbol.
\end{abstract}

\markright{Полимультипликативные отображения группы петель}

\footnotetext[0]{The study has been funded within the framework of the HSE University Basic Research Program
.}

\section{Introduction}\label{s1}


In 1994, C. Contou-Carrère, in his work \cite{CC_{1}}, introduced for any associative commutative ring $A$ a functorial pairing
$$A((t))^{*} \times A((t))^{*} \longrightarrow A^{*},$$
which is bimultiplicative, antisymmetric, and continuous.

If $A$ is a field, then it coincides with the tame symbol. This pairing is now called the Contou-Carrère symbol, and we will denote it by $CC_{1}$. Later, Deligne, in his work \cite{Deligne}, provided an explicit formula for this symbol. To formulate it, let us recall a fact about invertible elements in Laurent series over a ring. For simplicity, we will assume that $A$ does not decompose in a non-trivial way into a product of two rings, which is equivalent to the connectedness of the topological space $\Spec \ A$. In this case, an element $f = \sum \limits_{i\geq n} a_{i}t^{i} \in A((t))$ is invertible if and only if there exist $v(f) \in \Z$ such that $a_{v(f)}$ is invertible in $A$ and $a_{i}$ is nilpotent for every $i<v(f)$. The number $v(f)$ is an analogue of the of the standard valuation for Laurent series over a field. Now suppose $f, g\in A((t))^{*}$ and $A$ is a $\Q$-algebra, such that $v(f) = 0$ and the constant term of $f$ is $1$.  Then the following formula holds
\begin{equation}
\label{Deligne}
CC_{1}(f, g) = \exp (\Res(\log(f) \frac{dg}{g})).
\end{equation}

We denote by $\Res$ the residue of a differential $1$-form. Moreover, the following formulas hold: $CC_{1}(t, t) =-1$, $CC_{1}(a, g) = a^{v(g)}.$ This completely determines the Contou-Carrère symbol.

The main motivation for its study is that it satisfies a reciprocity law. Namely, suppose that $C$ is a smooth projective curve over an algebraically closed field $F$, $\nobreak{f,g\in (F(C) \otimes_{F} A)^{*},}$ where $A$ is a local Artinian $F$-algebra. Choose for every closed point of the curve a local paramerer at this point in the ring $\mathcal{O}_{C, x}.$ This determines the map $(F(C) \otimes_{F} A)^{*} \to A((t))^{*}$. Let us denote $f_{x}, g_{x}$ the images of $f, g$ under this map respectively. Then the following equality holds: 
\begin{equation}
\prod\limits_{x \in C} CC_{1}(f_{x}, g_{x}) = 1.
\end{equation}

This reciprocity law was proved by A. Beilinson, S. Bloch and H. Esnault in \cite{BBE}, and by G. Anderson, F. Pablos Romo in \cite{AndersonRomo}. In their works, the Contou-Carrère symbol is interpreted using a certain central extension of the group $A((t))^{*}$ by $A^{*}$. This approach can be considered as a development of D. Tate's approach to the residues of differential forms (\cite{Tate}) and E. Arbarello, C. De Concini, and V. Kac's approach to the tame symbol (\cite{ArbConKac}).

The multidimensional Contou-Carrère symbol was first defined by D. Osipov and X. Zhu in \cite{OsipovZhu1}, \cite{OsipovZhu2} for $n = 2,$ and later by S. Gorchinskiy and D. Osipov for an arbitrary $n$ in \cite{GorOsi15} (see also \cite{GorOsi20short}). Namely, let us denote $\mathcal{L}^{n}(A) = A((t_{1}))...((t_{n})).$ Then there is a map which depends functorially on a ring $A$:
$$CC_{n}\colon (\mathcal{L}^{n}(A)^{*})^{\times(n+1)} \longrightarrow A^{*}.$$

It is polymultiplicative and antisymmetric. Moreover, it satisfies a formula analogous to \eqref{Deligne}:
\begin{equation}
CC_{n}(f_{0}, f_{1},...,  f_{n}) = \exp (\Res(\log(f_{0}) \frac{df_{1}}{f_{1}} \wedge ... \wedge \frac{df_{n}}{f_{n}})),
\end{equation}
where $A$ is a $\Q$-algebra, $f_{i} \in \mathcal{L}^{n}(A)^{*}$ and there are some conditions on $f_{0}$ so that series $\log(f_{0})$ converges.

In \cite{GorOsi15} the multidimensional Contou-Carrère symbol $CC_{n}$ is described as a unique up to taking powers functorial polymultiplicative map which factors through Milnor's $K$-group. Explicitly, this means the following. Suppose that for every ring $A$ we are given a map which is functorial by a ring:
$$(\cdot, ... , \cdot)\colon (\mathcal{L}^{n}(A)^{*})^{\times(n+1)} \longrightarrow A^{*}.$$

Moreover, suppose it satisfies the following \textit{Steinberg relations}: if $f_{0}, ... , f_{n} \in \mathcal{L}^{n}(A)^{*}$ such that $f_{i+1} = 1 - f_{i}$ for some $i$, then $(f_{0}, ... , f_{n}) = 1$. Then for some integer $m \in \Z$ we have $(\cdot, ... , \cdot) = CC_{n}^{m}.$

In this work, we provide an alternative characterization of the Contou-Carrère symbol. Namely, iterated Laurent series have a natural topology, and so it makes sense to consider continious automorphisms of $\mathcal{L}^{n}(A)$. In \cite{GorOsi16} it was shown that multidimensional Contou-Carrère symbol is invariant under continious automorphisms. Namely, for any continious automorphism $\varphi \colon \mathcal{L}^{n}(A) \to \mathcal{L}^{n}(A)$ and a collection of elements $f_{0}, ... , f_{n} \in \mathcal{L}^{n}(A)^{*}$ we have the following equality:
\begin{equation}
\label{invariant}
CC_{n}(\varphi(f_{0}), ... , \varphi(f_{n})) = CC_{n}(f_{0}, ... , f_{n}).  
\end{equation}

This is generalization of the well known property of the tame symbol. We will prove that, up to a sign, the converse is true as well. Namely, the main result of the paper is the following theorem.
\begin{theorem}
Let $R$ be a torsion-free ring(i.e. additive group of $R$ has no $\Z$-torsion). Suppose also that for every $R$-algebra $A$ we are given a polymultiplicative map which is functorial by a ring A:
$$(\cdot, ... , \cdot) \colon (\mathcal{L}^{n}(A)^{*})^{\times(n+1)} \longrightarrow A^{*}.$$

Suppose also that for every continious automorphism ${\varphi \colon \mathcal{L}^{n}(A) \to \mathcal{L}^{n}(A)}$ and any collection $f_{0}, ... , f_{n} \in \mathcal{L}^{n}(A)^{*}$ the following holds:
$$(\varphi(f_{0}), ... , \varphi(f_{n})) = (f_{0}, ... , f_{n}). $$

Then there exists an element $m \in \underline{\Z}(R)$ and polylinear functorial map $$\omega \colon (\underline{\Z}_{R}^{n})^{\times(n+1)} \longrightarrow \G_{m, R}$$
such that for any ring $A$ and a collection $f_{0}, ... , f_{n} \in \mathcal{L}^{n}(A)^{*}$ the following equality holds:
\begin{equation}
\label{maineq}
(f_{0}, ... , f_{n}) = \omega(v(f_{0}), ... , v(f_{n})) CC_{n}(f_{0}, ... , f_{n})^{m}.
\end{equation}
\end{theorem}

We denote by $\underline{\Z}_{R}$ a sheaf associated to the constant presheaf. For example, in the case $R = \Z$ one can see that $\omega$ takes values in $\pm 1.$ So, in this case, we obtain that $(\cdot, ... , \cdot)$ coincides, up to a sign, with an integer power of $CC_{n}$. An element $v(f)\in \underline{\Z}^{n}(A)$ is the multidimensional analogue of the valuation from the one-dimensional case. We will recall these definitions in \S\ref{maindef}.

Actually we will prove the theorem under weaker assumptions. Namely, it is sufficient to require that, when restricted to $R\otimes \Q$-algebras, the map $(\cdot, ... , \cdot)$ is invariant under automorphisms $\varphi$, for which $\varphi(t_{i}) \in \mathcal{L}^{n}(\Q).$

The main technical ingredient of the proof is the following additive analogue, which is of independent interest. Namely, let $R$ be a ring which has sufficiently many invertible elements (see theorem \ref{theorem1} for more precise formulation) and let
$$<\cdot, ... , \cdot>\colon (\mathcal{L}^{n}(R))^{\times (n+1)} \longrightarrow R$$
be an $R$-polylinear map which is continious in each argument. Suppose further that it is  invariant under continious automorphisms $\varphi$ (in a sense analogous to \eqref{invariant}) for which $\varphi(t_{i}) \in \mathcal{L}^{n}(\Z)$ and sufficiently many automorphisms of the form $t_{i} \mapsto a_{i}t_{i}.$ Then there exist an element $e\in R$ such that the following equality holds:
$$<f_{0}, ..., f_{n}> = e \Res(f_{0} d f_1 \wedge...\wedge d f_{n}).$$

In the case $n = 1$ or $n = 2$ this statement has been proven in \cite{OsipovZhu2} (see lemmas 5.7, 5.8) to derive explicit formulas for $CC_{1}$ and $CC_{2}$ respectively.

This paper is organized as follows. In \S\ref{maindef} we recall the necessary definitions and facts about iterated Laurent series. Then in \S\ref{Laddth} we prove theorem \ref{theorem1} which says that polylinear map from iterated Laurent series, that is invariant under continious automorphisms, is the residue of a differential form up to a constant. In \S\ref{functorialprop} we prove some lemmas on iterated Laurent series which will be of use to prove the main theorem. Finally, in \S\ref{proofs} we prove the main theorem \ref{mainthm}, which is the main result of this work.

The author would like to express sincere gratitude to his supervisor, Denis Osipov, for suggesting the problem and for his attentive guidance throughout this work.

\section{Main definitions}\label{maindef}
In this chapter we only recall some facts about iterated Laurent series. Proofs can be found in \cite{GorOsi15}, \cite{GorOsi16}, \cite{GorOsi20}.

\subsection{Iterated Laurent series}

By a ring we always mean an associative commutative ring with a unit. For a ring $A$ we will denote by $\Nil(A)$ the ideal of nilpotent elements. By a group functor we mean a covariant functor from the category of rings to the category of groups. For example, one can define the functor $\Nil,$ which sends a ring $A$ to the abelian group $\Nil(A).$ We will identify an affine scheme with the functor it represents. 

Let $X$ be a functor from the category of rings to the category of sets. Recall that it is called an affine ind-scheme if it can be represented as a colimit of affine schemes $X_{i}$, $i\in I$, where $I$ is a \textit{directed} partially ordered set. More explicitely, $X = \varinjlim\limits_{i \in I} X_{i}$. One can see that it is the same as an ind-object in the category of affine schemes. Obvious generalizations can be made if we replace the category of rings with the category of $A$-algebras where $A$ is some ring. If X is a functor from $A$-algebras and $B$ is any $A$-algebra we denote by $X_{B}$ the restriction of $X$ to the category of $B$-algebras. For an ind-scheme $X$ over a ring $R$ we denote by $\mathcal{O}(X)$ the ring of regular functions on $X$, which can be defined as the set of natural transformations from $X$ to $\mathbb{A}^{1}_{R}.$

Let $G_{1}, G_{2}$ be commutative group functors over a ring $R$. We will denote $\Homgr(G_{1}, G_{2})$ a group functor which sends a $R$-algebra $A$ to the set of group morphisms grom $(G_{1})_{A}$ to $(G_{2})_{A}$, equipped with componentwise product. Also, we will denote by $G_{1}\otimes G_{2}$ a group functor that sends an $R$-algebra $A$ to the group $G_{1}(A) \otimes_{\Z} G_{2}(A)$. One can verify that there is an adjunction
$$\Homgr(G_{1}\otimes G_{2}, G_{3}) = \Homgr(G_{1}, \Homgr(G_{2}, G_{3})).$$ 

Let $A$ be a ring. By a ring of Laurent series we mean a set of formal power series of the form $\sum\limits_{i \geq i_{0}} a_{i}t^{i}$. Inductively, one defines a ring of \textit{iterated} Laurent series over $A$ that depends on the local parameters $t_{1}, ... , t_{n}$. Namely, we define $A((t_{1}))...((t_{n})) = (A((t_{1}))...((t_{n-1})))((t_{n}))$ equipped with the obvious operations of summation and multiplication. 

For a vector $l = (l_{1}, ... , l_{n}) \in \Z^{n}$ we denote $t^{l} = t_{1}^{l_{1}}... t_{n}^{l_{n}}.$ Thus, any element $f \in A((t_{1}))...((t_{n}))$ can be uniquely represented as an infinite sum of monomials: $f = \sum a_{l}t^{l}$. For such an element $f$ we denote by supp$(f) \subset \Z^{n}$ the set of vectors $l$, such that $a_{l} \neq 0.$ 

One can define a functor from the category of rings to the category of rings that sends a ring $A$ to the ring $A((t_{1}))...((t_{n}))$, equipped with the obvious operations of summation and multiplication. Moreover, for a group functor $G$, we define the functor $\mathcal{L}(G)$ that sends a ring $A$ to the group $G(A((t))).$ It is easy to see that $\mathcal{L}^{n}\G_{a}(A) = A((t_{1}))...((t_{n})).$ We also denote by $\mathcal{L}^{n}(A)$ the ring $A((t_{1}))...((t_{n})).$

In \cite{GorOsi15}, it was shown that $\mathcal{L}^{n}\G_{a}$ is an affine ind-scheme. Let us outline the main constructions. Introduce the set:
$$\Lambda_{n} = \{(\lambda_{1}, ..., \lambda_{n}) \mid \lambda_{p} \colon \Z^{n-p} \to \Z, 1\leq p \leq n \}.$$ 

Here by $\Z^{0}$ we mean one-point set. Thus, $\lambda_{n}$ can be identified with an element from $\Z$.

Next, for $\lambda \in \Lambda_{n}$, define a set $\mathbb{Z}^{n}_{\lambda}$ as follows: 
$$ \mathbb{Z}^{n}_{\lambda} = \{(l_{1}, ..., l_{n}) \in \Z^{n} \mid l_{n} \geq \lambda_{n}, l_{n-1} \geq \lambda_{n-1}(l_{n}), ... , l_{1} \geq \lambda_{1}(l_{2}, ... , l_{n}) \}.$$

Finally, for $\lambda \in \Lambda_{n}$, define a subfunctor $\mathcal{L}^{n}(\G_{a})_{\lambda} \subset \mathcal{L}^{n}(\mathbb{G}_{a})$ as follows:
$$\mathcal{L}^{n}(\mathbb{G}_{a})_{\lambda}(A) = \left\{ f\in \Ladd(A) \mid \mathrm{supp}(f) \subset \Z^{n}_{\lambda} \right\}.$$

First, one can prove that $\mathcal{L}^{n}(\G_{a})(A) = \bigcup\limits_{\lambda \in \Lambda_{n}}\mathcal{L}^{n}(\G_{a})_{\lambda}(A).$ Second, the functor $\mathcal{L}^{n}(\G_{a})_{\lambda}$ is represented by
$$\mathcal{L}^{n}(\G_{a})_{\lambda} = \Spec \ \Z[t_{l}]_{l\in \Z^{n}_{\lambda}}.$$

We define a topological group structure on the set $\mathcal{L}^{n}(\G_{a})(A)$. To do this, it suffices to specify a base of neighborhoods of zero. If $n = 1,$ the topology is defined using the standard choice of neighborhoods of zero, $t^{m} A[[t]], m \in \Z.$ Next, if $n>1$, we define the topology on $\mathcal{L}^{n}(A)$ using a base of neighborhoods of zero consisting of $A$-submodules of the form:
$$ U_{m, \{V_{j}\}_{j\leq m}} = \bigoplus_{j<m}t_{n}^{j}V_{j}\oplus t_{n}^{m}A((t_{1}))...((t_{n-1}))[[t_{n}]], $$
where $V_{j}$ is a $A$-module which belongs to the set of the open neighbourhood base at zero of $A((t_{1}))...((t_{n-1})).$

In the work \cite{GorOsi15} (see Lemma 3.3), it was shown that the group $\mathcal{L}^{n}(\G_{a})(A)$ is separated and complete. More specifically, the underlying topological space of $\mathcal{L}^{n}(\G_{a})(A)$ is Hausdorff and every Cauchy sequence has a limit.  Moreover, using lemma 3.3(b), one can show that a countable sequence ${f_{k} = \sum a_{k,l} t^{l}, k \in \NN}$ converges to $f = \sum a_{l} t^{l}$ if and only if there exist an element $\lambda\in \Lambda_{n}$, such that $f_{k}, f \in \mathcal{L}^{n}(\G_{a})_{\lambda}$ and for every $l\in \Z^{n}$ the coefficient $a_{k, l}$ is equal to $a_{l}$ for $k$ sufficiently large.

We also note that for $n>1$ the set $\mathcal{L}^{n}(A)$ is not a topological ring. Specifically, multiplication fails to be continuous with respect to the topology introduced above. On the other hand, multiplication by an element $f \in \mathcal{L}^{n}(A)$ is always continious. Moreover, multiplication in $\mathcal{L}^{n}(A)$ is \textit{sequentially} continuous. This means that if the sequences $f_{k}, g_{k},$ $k\in \NN$ converge to $f, g$ respectively, then the sequence $f_{k}g_{k}$ converges to $fg$ in $\mathcal{L}^{n}(A).$  

We say that an element $f \in \mathcal{L}^{n}(A)$ is \textit{topologically nilpotent} if the sequence $f^{i}, i\in \NN,$ converges to zero. In \cite{GorOsi15} the explicit description of topologically nilpotent elements is given. 

In order to outline this description, we define a \textit{lexicographical} order on the set $\Z^{n}$. We say that $(l_{1}, ... , l_{n}) \leq (l_{1}', ... , l_{n}'),$ if $l_{n} < l_{n}'$, or $l_{n} = l_{n}'$ and $(l_{1}, ... , l_{n-1}) \leq (l_{1}', ... , l_{n-1}').$

Now, we define $\Ladd^{\#}(A)$ to be the set of topologically nilpotent elements of the ring $\mathcal{L}^{n}(A)$. Then the following equality holds 
$$\Ladd^{\#}(A) = \left\{ \sum a_{l}t^{l} \in \mathcal{L}^{n}(A) \biggm| \sum\limits_{l\leq 0}a_{l} t^{l} \in \Nil(\mathcal{L}^{n}(A))  \right\}.$$

From this description, one can derive that the functor sending a ring $A$ to the set of topologically nilpotent elements, is an affine ind-scheme. Explicitly, for every pair $\lambda \in \Lambda_{n}$ and $k \in \mathbb{N}$ define the set: 
$$\Ladd^{\#}_{k , \lambda}(A) = \left\{\sum\limits_{l \in \Z^{n}_{\lambda}} a_{l}t^{l} \biggm| \Bigl(\sum\limits_{l\leq 0}a_{l} t^{l}\Bigr)^{k} = 0  \right\}.  $$

Then, obviously, $\Ladd^{\#}(A) = \bigcup\limits_{k, \lambda} \Ladd^{\#}_{k , \lambda}(A).$ Moreover, the subfunctor $\Ladd^{\#}_{k , \lambda}$ is the closed affine subscheme of $\Ladd_{\lambda}$ defined by a \textit{homogeneous} ideal $I_{k, \lambda}$, where the grading of each variable $t_{l}$ is $1.$

Now, we describe the invertible elements $\Lmult(A)$ in the ring $\mathcal{L}^{n}(A).$ From the fact that the group $\Ladd(A)$ is complete, it is easy to see that every element of the form $1 + f,$ where $f$ is topologically nilpotent, lies in $\Lmult(A).$

Define a group functor $\underline{\Z}$ that sends a ring $A$ to the group of continious maps from $\Spec \: A$ to $\Z,$ where we assume that the topology on $\Z$ is discrete. It is easy to see that $\underline{\Z}$ is the sheaf in the Zariski topology associated to the constant presheaf. Using the quasicompactness of the topological space $\Spec \: A,$ one immediately sees that defining such a map $\Spec \: A \to \Z$ is equivalent to giving a decomposition $A = \prod\limits_{i = 1}^{m} A_{i}$ together with a collection $k_{1}, ... , k_{m} \in \Z.$ Now, for $u \in \underline{\Z}(A)$ and $a = (a_{1}, ... , a_{m}) \in A$, we define $a^{u} = (a_{1}^{k_{1}}, ... , a_{m}^{k_{m}})$. Next, for $\mathcal{L}^{n}(A) = A((t_{1}))...(t_{n}))$, $l = (l_{1}, ... , l_{n}) \in \underline{\Z}^{n}(A),$ we put 
$$t^{l} = t_{1}^{l_{1}} ...  t_{n}^{l_{n}},$$
where we identify $l_{i}$ with its image in $\underline{\Z}(A((t))).$

If n = 1, in the work \cite{CC2} it was shown that every element $f \in \Lmult(A)$ can be uniquely decomposed in the form $f = a f_{+} f_{-} t^{l}$, where $a \in \G_{m}(A) = A^{*}$, $f_{+} \in 1 + tA[[t]]$, $f_{-} \in 1 + t^{-1} \Nil(A)[t^{-1}]$, $l \in \underline{\Z}(A)$. It can be formulated as an isomorphism of certain groups. We will state this immediately for the multidimensional case, following the work \cite{GorOsi15}. Specifically, there exists an isomorphism of group functors
\begin{equation}
\label{decomposition}
\Lmult =\G_{m} \times  \V_{+} \times \V_{-} \times \underline{\Z}^{n},
\end{equation}
where
$$\V_{+}(A) = \left\{ 1 + \sum\limits_{l>0} a_{l} t^{l} \in \mathcal{L}^{n}(A) \biggm| a_{l} \in A \right\}, $$
$$\V_{-}(A) = \left\{ 1 + \sum\limits_{l<0} a_{l} t^{l} \in \mathcal{L}^{n}(A) \biggm| \sum\limits_{l<0} a_{l} t^{l} \in \Nil(\mathcal{L}^{n}(A)) \right\}.$$

Note, that the product of an element from $\V_{+}$ by an element from $\V_{-}$ is always of the form $1 + f,$ where $f$ is topologically nilpotent. Denote by $v$ the projection from $\Lmult$ to the summand $\underline{\Z}^{n}$ in the decomposition \eqref{decomposition}. From this, it immediately follows that an element $f \in \Lmult(A)$ can be represented in the form $f = a (1 + g) t^{v(f)}$, where $g$ is topologically nilpotent and $a \in A^{*}$.

We will denote by $\Lmult^{\#}$ the group functor that sends a ring to the subgroup of the invertible elements of the form $1 + f,$ where $f$ is topologically nilpotent. More explicitely, from the discussion above, we see that
$$\Lmult^{\#}(A) = \left\{ 1 + \sum a_{l}t^{l} \in \mathcal{L}^{n}(A) \biggm| \sum\limits_{l\leq 0}a_{l} t^{l} \in \Nil(\mathcal{L}^{n}(A))  \right\}.$$

Moreover, $\Lmult^{\#}$ is an affine ind-scheme. Indeed, 
$$\Lmult^{\#} = \bigcup\limits_{k, \lambda} \Lmult^{\#}_{k, \lambda},$$
where
$$\Lmult^{\#}_{k, \lambda}(A) = \left\{1 + \sum\limits_{l \in \Z^{n}_{\lambda}} a_{l}t^{l} \biggm| \Bigl( \sum\limits_{l\leq 0}a_{l} t^{l} \Bigr)^{k} = 0  \right\}. $$

\subsection{Differential forms and Residues} \label{diffforms}
Let $\Omega_{\mathcal{L}^{n}(A)/A}^{1}$ be the module of Kähler differentials. We define a module of continuous Kähler differentials as follows:

$$\Omega^{\rm{cont}, 1}_{\mathcal{L}^{n}(A)/A} = \Omega^{1}_{\mathcal{L}^{n}(A)/A}/Q,$$
where $Q$ is generated by elements of the form $df - \sum\limits_{i=1}^{n} \frac{\partial f}{\partial t_{i}}dt_{i}, f\in \mathcal{L}^{n}(A).$ Further, by $df$, we always mean its image in $\Omega^{\rm{cont}, 1}_{\mathcal{L}^{n}(A)/A}.$ One can prove the following proposition.

\begin{propos}
\label{kahler}
The module $\Omega^{\rm{cont}, 1}_{\mathcal{L}^{n}(A)/A}$ is a free $\mathcal{L}^{n}(A)$-module of rank $n$, with basis $dt_{1}, ... , dt_{n}$.
\end{propos}

Now, we define $\Omega^{\rm{cont}, k}_{\mathcal{L}^{n}(A)/A}$ as the $k$-th wedge power of the module $\Omega^{\rm{cont}, 1}_{\mathcal{L}^{n}(A)/A}.$ From proposition \ref{kahler} it follows that
$$\Omega^{\rm{cont}, k}_{\mathcal{L}^{n}(A)/A} = \bigoplus\limits_{i_{1}< ... < i_{k}} \mathcal{L}^{n}(A) dt_{i_{1}}\wedge ... \wedge dt_{i_{k}}.$$

\begin{definition}
Let $\omega \in \Omega^{\rm{cont}, n}_{\mathcal{L}^{n}(A)/A}$. Write $\omega  = f dt_{1} \wedge ... \wedge dt_{n},$ where $f = \sum a_{l}t^{l}\in \mathcal{L}^{n}(A).$ Then we put $\Res(\omega) = a_{-1, -1, ... , -1}.$
\end{definition}

We now define an $A$-multilinear form on the additive group , depending on $n+1$ variables,
$$\Res \colon \Ladd(A)^{\times (n+1)} \longrightarrow A$$
by the formula:
$$\Res(f_{0}, ... , f_{n}) = \Res(f_{0} df_{1} \wedge ... \wedge df_{n}). $$

It is easy to verify that this form is continuous. Additionally, this form is invariant under continuous automorphisms. The proof can be found in \cite[Corollary 5.4]{GorOsiCont}.

\subsection{Exponent and logarithm}
Later, we will need the definition of the exponent and logarithm of iterated Laurent series.

\begin{propos}[\cite{GorOsi15}, Proposition 4.9]
On the category of $\Q$-algebras, there exists an isomorphism of group functors:
    $$ \exp \colon \Ladd_{\Q}^{\#} \longrightarrow \Lmult_{\Q}^{\#},$$
defined by the standard series:
$$ \exp(f) = \sum\limits_{k = 0}^{\infty} \frac{f^{k}}{k!}. $$

Its inverse is given by the series:
$$\log(1+f) = \sum\limits_{k=0}^{\infty} (-1)^{k+1}\frac{f^{k}}{k}.$$
\end{propos}

We note, however, that for $n>1$ these maps are \textit{not} continuous in the restriction topology. For example, let $n=2$ and $A = \Q[x_{1}, x_{2}, ...],$ where $x_{k}^{k} = 0$. Then it is easy to verify that the sequence $f_{k} = x_{k}t_{2}^{-1}t_{1}^{k}$ converges to zero, but the sequence $\exp(f_{k})$ does not have a limit.

\subsection{Higher-dimensional Contou-Carrère symbol}

In the work \cite{GorOsi15}, a higher-dimensional Contou-Carrère symbol was defined using a boundary map in $K$-theory. Namely, the higher-dimensional Contou-Carrère symbol is the polymultiplicative functor morphism:
$$ CC_{n} \colon (\Lmult)^{\times(n+1)} \longrightarrow \G_{m}.$$

In order to describe it explicitely, let us define a sign map.

\begin{propos}[\cite{GorOsi15}, Proposition 8.15]
There exists a unique polylinear, antisymmetric morphism of abelian groups:
$$ \sgn \colon \Z^{n} \times ... \times \Z^{n} \longrightarrow \Z / 2\Z,$$
depending on $n+1$ variables, such that for any $l, s_{2},..., s_{n} \in \Z^{n}$ the following equality holds:
$$\sgn(l, l, s_{2}, ... , s_{n}) = \det(l, s_{2}, ... , s_{n}).  $$

\end{propos}

\begin{remark}
The question of an explicit formula for this sign arose in the work of A.N. Parshin on the local class field theory(see \cite{Par84}). The answer was provided by S.V. Vostokov and I.B. Fesenko (see \cite{Par90}).Specifically, the following formula holds:
\begin{equation}
\sgn(l_{1}, ... , l_{n+1}) =\sum\limits_{1 \leq i < j \leq n} \det(l_{1}, ... , \hat{l_{i}}, ..., \hat{l_{j}}, ... , l_{n+1}, l_{i}\cdot l_{j}),
\end{equation}
where $l_{i}\cdot l_{j}$ is the componentwise multiplication of elements in $\Z^{n}$, and the right-hand side of the equation is taken modulo $2.$ Another formula is due to A.G. Khovanskii (see \cite[Theorem 1]{Kho06}, see also \cite[Proposition 11]{Osipov2003} for the case $n = 2$):
\begin{equation}
\begin{aligned}
\sgn(l_{1}, ... , l_{n+1}) = 1 + \sum\limits_{i = 1}^{n+1} \det(l_{1}, ... , \hat{l_{i}}, ... , l_{n+1}) + \\ + \prod\limits_{i = 1}^{n+1} \left(1 + \det(l_{1}, ... , \hat{l_{i}}, ... , l_{n+1}) \right).
\end{aligned}
\end{equation}
\end{remark}

The map $CC_{n}$ is completely described by the following theorem:

\begin{theorem}[\cite{GorOsi15}, Theorem 8.17, Proposition 8.22]
The functorial morphism $CC_{n}$ satisfies the following properties:
\begin{itemize}
\item $CC_{n}$ is polymultiplicative.
\item If $A$ is an arbitrary ring, $a\in \G_{m}(A),$ $f_{1}, ... , f_{n} \in \Lmult(A),$ then
\begin{equation} \label{CCeq1}
CC_{n}(a, f_{1}, ... , f_{n}) =  a^{\det(v(f_{1}), ... , v(f_{n}))}.
\end{equation}
\item If $A$ is a $\Q$-algebra, $f_{0} \in \Lmult(A)^{\#}, f_{1}, ... , f_{n} \in \Lmult(A),$ then 
\begin{equation}\label{CCeq2}
CC_{n}(f_{0}, ... , f_{n}) = \exp \Res (\log f_{0} \frac{df_{1}}{f_{1}} \wedge ... \wedge \frac{df_{n}}{f_{n}}).
\end{equation}
\item For any collection of vectors $l_{0}, ... , l_{n} \in \Z^{n}$, the following holds:
\begin{equation}\label{CCeq3}
CC_{n}(t^{l_{0}}, t^{l_{1}}, ..., t^{l_{n}}) = (-1)^{\sgn (l_{0}, ... , l_{n})}.
\end{equation}
\end{itemize}
Moreover, these properties uniquely determine $CC_{n}$.
\end{theorem}

In what follows, by $CC_{n}$ we mean a morphism $(\Lmult)^{\times(n+1)}_{\Q} \longrightarrow \G_{m, \Q}$, that is polymultiplicative, antisymmetric, and satisfies \ref{CCeq1}, \ref{CCeq2}, \ref{CCeq3}.

\subsection{Continuous automorphisms} \label{contauto}
Now we describe continuous ring automorphisms of the $\mathcal{L}^{n}(A)$. The proofs can be found in \cite{GorOsiCont}. First, let $\varphi \colon \mathcal{L}^{n}(A) \to \mathcal{L}^{n}(A)$ is some continuous ring homomorphism. Then put
$$M(\varphi) = (v(\varphi(t_{1})), ... , v(\varphi(t_{n}))) \in Mat_{n\times n}(\underline{\Z}(A)).$$

The following formula holds:
\begin{equation}\label{eqauto}
v(\varphi(f)) = M(\varphi)v(\varphi(f))
\end{equation}

Moreover, for any two continuous ring homomorphisms $\varphi_{1}, \varphi_{2}$ the following holds: $M(\varphi_{1} \varphi_{2}) = M(\varphi_{1}) M(\varphi_{2}).$ There exists the following criterion:
\begin{theorem}[\cite{GorOsiCont}, Theorem 6.8]
\label{auto}
If $\varphi$ is a continuous ring automorphism of $\mathcal{L}^{n}(A),$ then the matrix $M(\varphi)$ is upper-triangular with units on the diagonal. Conversely, if $f_{1}, ... , f_{n} \in \mathcal{L}^{n}(A)^{*}$ is a collection such that the matrix $(v(f_{1}), ... , v(f_{n}))$ is upper-triangular with units on the diagonal, then there exists a unique continuous automorphism $\varphi$ such that $\varphi(t_{i}) = f_{i}$.
\end{theorem}

In what follows, we will often use this result without providing explicit references. Note that this theorem implies the existence of a group functor that assigns to a ring $A$ the group of continuous $A$-linear automorphisms of the ring $\mathcal{L}^{n}(A).$ We denote this group by $Aut(\mathcal{L}^{n}(A))$. Furthermore, if $A$ is an $\Q$-algebra, then $Aut(\mathcal{L}^{n}(\Q))$ embeds in $Aut(\mathcal{L}^{n}(A))$.

For a ring $R$ we denote by $B(n, R)$ the group of upper-triangular matrix with coefficients in the ring $R$ and with units on the diagonal.

From the theorem \ref{auto}, it follows that there exists a group homomorphism
$$M \colon Aut(\mathcal{L}^{n}(A)) \to B(n, \underline{\Z}(A)),$$
which sends a continious automorphism to the matrix $M(\varphi).$ This homomorhism has a group section that assigns to a matrix $U$ the automorphism mapping $t^{l}$ to $t^{U(l)},$ where $U(l)$ is the product of the matrix $U$ with the column vector $l$. We denote this automorphism by the same letter $U$.

\section{Polylinear maps to the additive group of Laurent series} \label{Laddth}
In this section we assume that the topology on a ring $R$ is discrete. Recall that a group $H$ has an infinite exponent if for every natural $n>0$ there exists an element $g \in H,$ such that $g^{n} \neq 1.$

Let us proceed to the proof of the following theorem:
\begin{theorem}
\label{theorem1}
Let $R$ be a ring and $H\subset R^{*}$ be a subgroup of infinite exponent. Let
$$<\cdot, ... , \cdot> \colon (\Ladd(R))^{\times(n+1)} \longrightarrow \G_{a}(R)$$
be an $R$-polylinear map that is continuous in each argument and invariant under automorphisms from $Aut(\mathcal{L}^{n}(\Z))$ and automorphisms of the form $t_{i} \mapsto a_{i}t_{i}$, where $a_{i} \in H.$ Then there exists an element $e \in R,$ such that for any $f_{0}, ... , f_{n}\in \Ladd(R)$ the following holds:
$$<f_{0}, ... , f_{n}> = e \Res(f_{0} df_{1} \wedge ... \wedge df_{n}).$$
\end{theorem}

It was proved for the cases $n = 1$ and $n = 2$ in the work \cite{OsipovZhu2} in the case of a field of characteristic $0$. In our case, the proof is the same, so we use the case $n = 1$ without proof.

Let $f_{0}, ..., f_{n}$ be a collection of monomials of the form $f_{i} = t_{1}^{\alpha_{i, 1}}...t_{n}^{\alpha_{i, n}}$ satisfying the condition:
\begin{equation}
    \label{s1}
    \sum \limits_{i = 0}^{n} \alpha_{i, j} = 0.  
\end{equation}

Informally, the sum of degrees in each $t_{j}$ is zero. For such a data, we denote by $D(f_{0}, ..., f_{n})$ the $n \times n$ matrix composed of the elements $\alpha_{i,j},$ where $1 \leq i \leq n$, $1\leq j \leq n.$

\begin{lemma}
\label{residue}
 Let $f_{i} = t_{1}^{\alpha_{i, 1}}...t_{n}^{\alpha_{i, n}}$ be an arbitrary collection of monomials. Then:
 \begin{enumerate}
     \item If the condition (\ref{s1}) is satisfied, then 
	$$\Res(f_{0} df_{1}\wedge...\wedge df_{n}) = \det(D(f_{0},..., f_{n})).$$
     \item Otherwise, the following holds: 
	$$\Res(f_{0} df_{1}\wedge...\wedge df_{n}) = 0.$$
     
 \end{enumerate}
\end{lemma}

\begin{proof}
Taking derivatives one easily sees:
$$ df_{i} = f_{i}\sum\limits_{j = 1}^{n} \alpha_{i,j} \frac{dt_{j}}{t_{j}}. $$

Then
$$f_{0}df_{1}\wedge... \wedge df_{n} = \frac{f_{0}...f_{n}}{t_{1}...t_{n}} \det(D(f_{0},..., f_{n})) dt_{1}\wedge... dt_{n}.$$

From this, both parts of the lemma follow.
\end{proof}

\begin{lemma}
\label{restriction}
Let $R$ be a ring and $H\subset R^{*}$ is a subgroup of infinite exponent. Let
$$<\cdot, ... , \cdot>\colon (\Ladd(R))^{\times(n+1)} \longrightarrow \G_{a}(R)$$
be an $R$-polylinear map that is invariant under automorphisms of the form $t_{i} \mapsto a_{i}t_{i}$, where $a_{i} \in H.$ Then, if $f_{0}, ..., f_{n}$ are monomials of the form $f_{i} = t_{1}^{\alpha_{i, 1}}...t_{n}^{\alpha_{i, n}}$ and coefficients $\alpha_{i,j}$ do not satisfy $(\ref{s1}),$ then 
$$<f_{0}, ..., f_{n}> = 0. $$
\end{lemma}
\begin{proof}
Indeed, suppose the sum of degrees in some of the variables, say in $t_{i}$, is equal to $N \neq 0$. Then we can construst an automorphism that maps $t_{i}$ to $t_{i}' = at_{i},$ where $a\in H,$ and $t_{j} \mapsto t_{j}$, for $j\neq i.$ Then, by the conditions of the lemma, the following holds:
$$<f_{0},..., f_{n}> = a^{N}<f_{0}, ... f_{n}>.$$

Again, by the conditions of the lemma, there exists $a\in H, $ such that $a^{N} \neq 1$. Thus, we have $<f_{0}, ... , f_{n}> = 0$.
\end{proof}
\begin{lemma}
\label{restriction2}
Let $R$ be a ring and $H\subset R^{*}$ be a subgroup of infinite exponent. Let
$$<\cdot, ... , \cdot>\colon (\Ladd(R))^{\times(n+1)} \longrightarrow \G_{a}(R)$$
be an $R$-polylinear map that is invariant under automorphisms from $Aut(\mathcal{L}^{n}(\Z))$ and automorphisms of the form $t_{i} \mapsto a_{i}t_{i}$, where $a_{i} \in H.$ Additionally, suppose that $f_{i} = t_{1}^{\alpha_{i, 1}}...t_{n}^{\alpha_{i, n}}$ satisfy condition (\ref{s1}). If all the $f_{i}$ does not depend on some variable $t_{j}$, i.e. $\alpha_{i,j} = 0$ for every $i$, then $$<f_{0}, ..., f_{n}> = 0.$$
\end{lemma}
\begin{proof}
Suppose that all $f_{j}$ do not depend on $t_{j}.$ Then
\begin{equation}
    \label{r1}
    <f_{0}t_{j}^{-1}, f_{1}, ..., f_{n}> = 0,
\end{equation}
by lemma \ref{restriction}, since the sum of the degrees in $t_{j}$ is equal to $-1$. Next, we define an automorphism that maps $t_{k}$ to $t_{k}'$, where 
\begin{align*}
t_{j}' &= t_{j}(1 + t_{j})^{-1}, \\
t_{k}' &= t_{k}, \ \ k \neq j.
\end{align*} 

Substituting into (\ref{r1}), we obtain:
$$<f_{0}, ..., f_{n}> + <f_{0}t_{j}^{-1}, ..., f_{n}> = 0 .$$

From this, it follows that $<f_{0}, ..., f_{n}> = 0$.
\end{proof}

\begin{lemma}
\label{maintech}
Let $R$ be a ring and
$$(\cdot, \cdot) \colon R((t_{1}))((t_{2})) \times R((t_{1}))((t_{2})) \longrightarrow R$$
continuous $R$-bilinear map.
\begin{enumerate}
\item Suppose that there exists $p\in \Z$ with the propery that if $\omega + \gamma \neq p$ then $(t_{1}t_{2}^{\omega}, t_{1}^{-1}t_{2}^{\gamma}) = 0$. Additionaly, suppose that $(\cdot, \cdot)$ is invariant under automorphisms from $Aut(\mathcal{L}^{2}(\Z))$ that map $t_{2}$ to $t_{2}.$ Then, for all $\omega, \gamma \in \Z$ the following equality holds:
$$(t_{1}t_{2}^{\omega}, t_{1}^{-1}t_{2}^{\gamma}) = (t_{1}, t_{1}^{-1}t_{2}^{\omega + \gamma}).$$
\item Suppose that there exists $p\in \Z,$ with the propery that if $\omega + \gamma \neq p$ then $(t_{1}^{\omega}t_{2}, t_{1}^{\gamma}t_{2}^{-1}) = 0$. Additionaly, suppose that $(\cdot, \cdot)$ is invariant under automorphisms from $Aut(\mathcal{L}^{2}(\Z))$ that map $t_{1}$ to $t_{1}$. Then, for all $\omega, \gamma \in \Z$ the following equality holds:
$$(t_{1}^{\omega}t_{2}, t_{1}^{\gamma}t_{2}^{-1}) = (t_{2}, t_{1}^{\omega + \gamma}t_{2}^{-1}).$$ 
\item Let $H\subset R^{*}$ be a subgroup of infinite exponent. Suppose that $(\cdot, \cdot)$ is invariant under continuous automorphisms from $Aut(\mathcal{L}^{2}(\Z))$ and automorphisms of the form $t_{i} \mapsto a_{i}t_{i}$ where $a_{i}\in H$ that map $t_{2}$ to $t_{2}$. Then, for all $\omega, \gamma \in \Z$, the following holds:
$$(t_{1}t_{2}^{\omega}, t_{1}^{-1}t_{2}^{\gamma}) = -(t_{1}^{-1}t_{2}^{\omega}, t_{1}t_{2}^{\gamma}).$$
\item Let $H\subset R^{*}$ be a subgroup of infinite exponent. Suppose that $(\cdot, \cdot)$  is invariant under continuous automorphisms from $Aut(\mathcal{L}^{2}(\Z))$ and automorphisms of the form $t_{i} \mapsto a_{i}t_{i}$ where $a_{i}\in H$ that map $t_{1}$ to $t_{1}$. Then, for all $\omega, \gamma \in \Z$, the following holds:
$$(t_{1}^{\omega}t_{2}, t_{1}^{\gamma}t_{2}^{-1}) = -(t_{1}^{\omega}t_{2}^{-1}, t_{1}^{\gamma}t_{2}).$$
\end{enumerate}
\end{lemma}
\begin{remark}
Parts 1 and 2 in the statement of the lemma \ref{maintech} differ only by the permutation of $t_{1}$ and $t_{2}$. We include them both because, in the ring $R((t_{1}))((t_{2}))$, these elements are not symmetric. A similar remark applies to parts 3 and 4.
\end{remark}
\begin{proof}
Let us prove the first part. From a proof it will be clear that the second part can be proved analogously.

Fix elements $\omega, \gamma \in \Z.$ If $\omega + \gamma \neq p,$ then both sides of the equality is zero by conditions of the lemma. Now, if $\omega = 0$ or $\gamma = 0,$ then the equality is obvious. In what follows, we suppose that $\gamma, \omega \neq 0$ and $\omega + \gamma = p.$ Let us consider separately the cases when $\omega > 0$ and $\omega < 0$. 

Suppose that $\omega > 0$. Then there exists an automorphism that maps $t_{1}$ to $t_{1}(1 + t_{2}^{\omega})$ and $t_{2}$ to $t_{2}$. Then, by conditions of the lemma, we have
$$0 = (t_{1}, t_{1}^{-1}t_{2}^{\gamma}) = (t_{1}(1 + t_{2}^{\omega}), t_{1}^{-1}(1 + t_{2}^{\omega})^{-1}t_{2}^{\gamma}). $$

Using the series expansion for $(1 + t_{2}^{\omega})^{-1}$, expanding the brackets, and applying the conditions of the lemma, we obrain that
$$ 0 =  (t_{1}t_{2}^{\omega}, t_{1}^{-1}t_{2}^{\gamma}) - (t_{1}, t_{1}^{-1}t_{2}^{\omega + \gamma}).$$

Next, suppose that $\omega < 0.$ Then there exists an automorphism that sends $t_{1}$ to $\linebreak{t_{1}(1 + t_{2}^{-\omega})}$ and $t_{2}$ to $t_{2}$. Again, by conditions of the lemma, we derive that
$$0 = (t_{1}t_{2}^{\omega}, t_{1}^{-1}t_{2}^{\gamma + \omega}) = (t_{1}(1 + t_{2}^{-\omega})t_{2}^{\omega}, t_{1}^{-1}(1 + t_{2}^{-\omega})^{-1}t_{2}^{\gamma + \omega}).$$

Analogous to the case $\omega > 0$, by expanding the brackets in the expression above, we obtain:
$$0 = -(t_{1}t_{2}^{\omega}, t_{1}^{-1} t_{2}^{\gamma}) + (t_{1}, t_{1}^{-1}t_{2}^{\gamma + \omega}). $$ 

Let us prove the third part. Again, from a proof it will be clear that the fourth part can be proved analogously. Consider a map
$$<\cdot, \cdot> \colon R((t_{1})) \times R((t_{1})) \longrightarrow R,$$
given by the formula $<f,g> = (f t_{2}^{\omega}, g t_{2}^{\gamma}).$ It is continuous, $R$-bilinear and invariant under continuous automorphisms from $Aut(\mathcal{L}^{2}(\Z))$ and automorphisms of the form $t_{1} \mapsto at_{1},$ where $a\in H.$ Then we apply the case $n = 1$ of theorem \ref{theorem1} to obtain $<t_{1}, t_{1}^{-1}> = \linebreak{- <t_{1}^{-1}, t_{1}>}$. From this, the lemma follows.
\end{proof}

\begin{lemma}
\label{permutations}
Let $R$ be a ring and $H\subset R^{*}$ be a subgroup of infinite exponent. Let
$$<\cdot, ... , \cdot>\colon (\Ladd(R))^{\times(n+1)} \longrightarrow \G_{a}(R)$$
be an $R$-polylinear map that is continuous in each argument and invariant under automorphisms from $Aut(\mathcal{L}^{n}(\Z))$ and automorphisms of the form $t_{i} \mapsto a_{i}t_{i}$, where $a_{i} \in H.$ Then if $\pi \in S_{n+1}$ is an arbitrary permutation, $y_{0} = t_{1}^{-1}...t_{n}^{-1}$ and $y_{i} = t_{i}$ for $1 \leq i \leq n,$ then the following equality holds:
$$ <y_{\pi(0)}, ..., y_{\pi(n)}> = \sgn(\pi) <y_{0}, ..., y_{n}> .$$
\end{lemma}

\begin{proof}
First, we prove that in the case $n = 2$, we can permute $y_{1}$ and $y_{2}$. More formally, we need to show the following equality:
$$<t_{1}^{-1}t_{2}^{-1}, t_{1}, t_{2}>=-<t_{1}^{-1}t_{2}^{-1}, t_{2}, t_{1}>.$$

We define $(\cdot, \cdot) \colon R((t_{1}))((t_{2}))^{\times 2} \to R$ by the formula $(f,g) = <g, f, t_{2}>$. From the first part of Lemma \ref{maintech} we obtain the following equality$$<t_{1}^{-1}t_{2}^{-1}, t_{1}, t_{2}> = (t_{1}, t_{1}^{-1}t_{2}^{-1}) = (t_{1}t_{2}^{-1}, t_{1}^{-1}) =  <t_{1}^{-1}, t_{1}t_{2}^{-1}, t_{2}>.$$

Using analogous considerations we obtain
$$<t_{1}^{-1}, t_{1}t_{2}^{-1}, t_{2}> = <t_{1}^{-1}, t_{2}^{-1}, t_{1} t_{2}> = - <t_{1}^{-1}, t_{2}, t_{1}t_{2}^{-1}> = -<t_{1}^{-1}t_{2}^{-1}, t_{2}, t_{1}>.$$

Now consider an arbitrary $n\geq 3$. Let us prove that we can permute $y_{1}$ and $y_{2}.$ Consider the following map 
$$(\cdot, \cdot, \cdot) \colon R((t_{1}))((t_{2}))^{\times 3} \longrightarrow R $$
given by the formula $(f,g,h) = (ft_{3}^{-1}...t_{n}^{-1}, g, h, t_{3}, ..., t_{n}).$ It is $R$-polylinear and continuous in each argument. Moreover, it is invariant under automorphisms from $Aut(\mathcal{L}^{2}(\Z))$ and automorphisms of the form $t_{i} \mapsto a_{i}t_{i}$ where $a_{i} \in H$. Thus, using what has already been proven, we obtain $(t_{1}^{-1} t_{2}^{-1}, t_{1}, t_{2}) = -(t_{1}^{-1} t_{2}^{-1}, t_{2}, t_{1})$. Now, using the definition of $(f,g,h)$, we see that we can permute $y_{1}$ and $y_{2}$ for every $n$. Then it is obvious that, by similar arguments, we can permute any $y_{i}, y_{j}$ for $i>j\geq 1.$

It suffices to show the equality
$$ <y_{0}, y_{1}, ..., y_{n}> = - <y_{n}, y_{1}, y_{2}, ..., y_{0}>. $$

There exists a continuous automorphism that maps $t_{i}$ to $t_{i}'$ where
\begin{align*}
t_{n}' &= t_{n}t_{1}^{-1}...t_{n-1}^{-1}, \\
 t_{i}' &= t_{i},  i<n. 
\end{align*}

Thus, by the conditions of the lemma, we obtain
$$ <y_{0}, y_{1}, ..., y_{n}> = <t_{n}^{-1}, y_{1}, ...,y_{n-1}, t_{1}^{-1}t_{2}^{-1}... t_{n-1}^{-1}t_{n}>.$$

Again, define a map $(\cdot, \cdot) \colon R((t_{n})) \times R((t_{n})) \longrightarrow R$ by the formula $\linebreak{(f,g) = <f, y_{1},..., y_{n-1}, t_{1}^{-1}...t_{n-1}^{-1}g>}.$ Then $(\cdot,\cdot)$ is continuous and $R$-bilinear. \linebreak{Moreover,} it is invariant under automorphisms from $Aut(\mathcal{L}(\Z))$ and automorphisms of the form $t_{n} \mapsto at_{n}$, where $a\in H$. Now, from the case $n = 1$ of Theorem \ref{theorem1}, it follows that any such pairing is antisymmetric. Thus, $(f,g) = - (g, f)$. Using the definition of $(\cdot, \cdot),$ we obtain the lemma.
\end{proof}

Now we are ready to prove Theorem \ref{theorem1}.

\begin{proof}
We proceed by induction on $n.$

From continuity, we see that it suffices to consider only monomials of the form $f_{i} = t_{1}^{\alpha_{i, 1}}...t_{n}^{\alpha_{i, n}}$. Moreover, by Lemmas \ref{residue} and \ref{restriction}, we can we can assume that condition $(\ref{s1})$ is satisfied. Otherwise, both sides of the equality are zero.

Thus, we have to show the following equality:
$$ <f_{0}, ..., f_{n}> = <t_{1}^{-1}...t_{n}^{-1}, t_{1}, ..., t_{n}> \det(D(f_{0},...,f_{n})).$$

If $\alpha_{i,n} = 0$ for every $i$, we see that both sides of the equality are zero by Lemmas \ref{restriction2} and \ref{residue}, and we are done. Otherwise, there exists $i$ such that $\alpha_{i,n} \neq 0$. Moreover, by condition \eqref{s1}, we can assume that $\alpha_{i,n}>0$. Indeed, for any permutation $\pi \in S_{n+1}$ consider a new map $<f_{0}, ... , f_{n}>_{\pi} = <f_{\pi(1)}, ..., f_{\pi(n)}>$. Furthermore, if we prove the desired statement for the new map, then for the original map it will follow using Lemma \ref{permutations}. We will consider two cases:

\begin{enumerate}
\item Suppose that $\alpha_{n,n} = 1.$ There exists a continuous automorphism that maps $t_{i}$ to $t_{i}'$, where
\begin{align*} 
t_{n}' &= t_{n} t_{1}^{-\alpha_{n, 1}}...t_{n-1}^{-\alpha_{n, n-1}}, \\
t_{k}' &= t_{k}, \ k<n.
\end{align*}
Thus, we obtain $<f_{0}, ..., f_{n}> = <f_{0}',..., f_{n-1}', t_{n}>$, where $f_{i}'$ is the image under automorphism considered above. Moreover, one easily sees that
$$\det D(f_{0},...,f_{n}) = \det D(f_{0}', ...,f_{n-1}', t_{n}).$$ 

Hence, if $\alpha_{n,n} = 1$, we can assume that $f_{n} = t_{n}.$

Next, for $i\leq n-1$, write $f_{i} = f_{i}'t_{n}^{\beta_{i}}$, where $f_{i}'$ does not depend on $t_{n}$, and $f_{n} = t_{n}$. We define an $R$-polylinear map $(R((t_{1}))...(t_{n-1})))^{\times n} \to R$ by the formula:
$$(g_{0}, g_{1}, ..., g_{n-1}) = <g_{0}t_{n}^{\beta_{0}},..., g_{n-1}t_{n}^{\beta_{n-1}}, t_{n}>.$$

It is continuous and invariant under continuous automorphisms as in the conditions of the theorem. By induction, we obtain
$$<f_{0}, ..., f_{n}> = (f_{0}',..., f_{n-1}') = \det D(f_{0}, ..., f_{n})(t_{1}^{-1}... t_{n-1}^{-1}, t_{1}, ..., t_{n-1}).$$

By definition 
$$(t_{1}^{-1}... t_{n-1}^{-1}, t_{1}, ..., t_{n-1}) = <t_{1}^{-1}... t_{n-1}^{-1}t_{n}^{\beta_{0}}, t_{1}t_{n}^{\beta_{1}},..., t_{n-1}t_{n}^{\beta_{n-1}} ,t_{n}>.$$

Now it suffices to show that expression on the right is equal to $<t_{1}^{-1}...t_{n}^{-1}, t_{1}, ..., t_{n}>$. We define $\{ \cdot, \cdot \} \colon R((t_{1}))((t_{n}))^{\times 2} \to R$ by the formula 
$$\{ f, g \} = <f t_{2}^{-1}...t_{n-1}^{-1}, g, t_{2} t_{n}^{\beta_{2}}, ... ,t_{n-1} t_{n}^{\beta_{n-1}}, t_{n}>.$$

One can easily see that $\{ \cdot, \cdot \}$ satisfies the conditions of the first part of Lemma \ref{maintech}. From this, we obtain $\{t_{1}^{-1} t_{n}^{\beta_{0} + \beta_{1}}, t_{1} \} = \{t_{1}^{-1} t_{n}^{\beta_{0}}, t_{1} t_{n}^{\beta_{1}} \}.$ Using the definition of $\{ \cdot, \cdot \}$, we can assume that $\beta_{1} = 0$. Proceeding like this, we can assume that $\beta_{i} = 0$ for $i\geq 1.$ But, $\beta_{0} + ... \beta_{n} = 0,$ so $\beta_{0} = 0$, and we are done.

\item Now consider the case $\alpha_{n, n} > 1$. As before, we write ${f_{i} = f_{i}'t_{n}^{\beta_i}}$, where $\beta_{i} = \alpha_{i, n}$.  
Now let $\alpha_{n, n} > 1$. As before, we write $f_i = f_i' t_n^{\beta_i}$, where $\beta_i = \alpha_{i, n}$. Furthermore, for $0\leq i \leq n-1$ we put
 $$\gamma_{i} = <f_{0}'t_{n}^{\beta_{0}}, ..., f_{i}'t_{n}^{\beta_{i} + \beta_{n} - 1},..., f_{n}' t_{n}>,$$
 $$\gamma_{n} = <f_{0}'t_{n}^{\beta_{0}}, ..., f_{i}'t_{n}^{\beta_{i}},..., f_{n}' t_{n}^{\beta_{n}}>. $$ 
 
Our goal is to compute $\gamma_{n}$. Observe that
 \begin{equation}
 \label{r5}
 <f_{0}'t_{n}^{\beta_{0}}, ..., f_{i}'t_{n}^{\beta_{i}},..., f_{n}' t_{n}> = 0.
 \end{equation}

There exists a continuous automorphism that maps $t_{k}$ to $t_{k}'$ where
 $$t_{n}' = t_{n} + t_{n}^{\beta_{n}}, $$
$$t_{k}' = t_{k}, \ \  k\leq n-1. $$

Then, by substituting $t_{i}'$ for $t_{i}$ in \eqref{r5}, we obtain the equality: 
$$ \sum\limits_{i = 0}^{n-1} \beta_{i} \gamma_{i} + \gamma_{n} = 0. $$
 
Now, $\gamma_{i}$ for $i \leq n-1$ has already been computed above. Then, everything follows from the expansion of the determinant of the matrix along the last column.
 \end{enumerate}
\end{proof}

\section{Functorial properties of Laurent series}\label{functorialprop}
It is easy to see that a functor $\underline{\Z}_{R}$, introduced in \S\ref{maindef}, is an ind-affine scheme. Namely, we identify representable functor $\bigsqcup\limits_{-n \leq i \leq n} \Spec  \ R$ with a subfunctor of $\underline{\Z}_{R},$ which assigns to an $R$-algebra $A$ the set of continuous maps from $\Spec \ A$ to $\Z$ such that the image of any element lies in the interval from $-n$ to $n$. Then, using quasi-compactness of $\Spec \ A$, we see that
$$ \underline{\Z}_{R} = \varinjlim\limits_{n} \bigsqcup\limits_{-n \leq i \leq n} \Spec \ R.$$

\begin{lemma}
\label{morphismstoZ}
Let $G$ be an affine group scheme over a ring $R$. Then any group morphism $G \to \underline{\Z}_{R}$ is trivial. 
\end{lemma}
\begin{proof}
Indeed, by the Yoneda lemma, we see that any such morhphism must factor through $\bigsqcup\limits_{-n \leq i \leq n} \Spec  \ R$. Now, it is obvious that any such morphism is trivial.
\end{proof}

The following proposition is well-known, so we omit a proof.
\begin{propos}
\label{Gm}
Any group morphism $\G_{m, R} \to \G_{m, R}$ is of the form ${a \mapsto a^{l}}$ for some $l \in \underline{\Z}(R).$
\end{propos}

\begin{cons}
\label{trivial}
Let $R$ be a ring. Then the following holds
\begin{enumerate}
\item Any group morphism $\G_{m, R} \to \underline{\Z}_{R}$ is trivial.
\item Any bimultiplicative map $\G_{m,R} \times \G_{m,R} \to \G_{m,R}$ is trivial.
\item Any group morphism $\Lmult^{\#}_{R} \to \underline{\Z}_{R}$ is trivial.
\item Any bimultiplicative map $\Lmult^{\#}_{R} \times \G_{m, R} \to \G_{m, R}$ is trivial.
\end{enumerate}
\end{cons}

\begin{proof}
The first part follows immediately from Lemma \ref{morphismstoZ}. The second part follows from the first part and Proposition \ref{Gm}, since such a bilinear map induces a morphism: $\G_{m,R} \to \Homgr(\G_{m, R}, \G_{m,R}) = \underline{\Z}_{R}$. The fourth part follows from the third part in the same way. It remains to prove the third part.

It suffices to prove that for any $k \in \NN$ and $\lambda \in \Lambda_{n}$, the induced morphism $\Lmult^{\#}_{R, \lambda, k} \to \underline{\Z}_{R}$ maps everything to zero. Put $A = \O(\Lmult^{\#}_{R, \lambda, k})$. By the Yoneda lemma, it follows that giving a morphism from $\Spec \ A$ to $\underline{\Z}_{R}$ is the same as giving an element $ x \in \underline{\Z}(A)$. However, $\underline{\Z}(A) = \underline{\Z}(A^{\mathrm{red}}),$ where $A^{\mathrm{red}} = A/\Nil(A).$ From the explicit description of topologically nilpotent elements it immediately follows that $A^{\mathrm{red}}$ is simply a polynomial algebra over $R$.

Now, as in Lemma \ref{morphismstoZ}, $x$ comes from $\bigsqcup\limits_{-n \leq i \leq n} \Spec  \ R.$ That is, it is determined by a set of idempotents $e_{i} \in A^{\mathrm{red}},$ where $ -n \leq i \leq n.$ However, since $A^{\mathrm{red}}$ is a polynomial algebra, such idempotents come from $R.$ Therefore, since the original morphism $\Lmult^{\#}_{R} \to \underline{\Z}_{R}$ was a group morphism, the statement follows.
\end{proof}

\begin{lemma}
\label{nilppolylinear}
Let $R$ be a $\Q$-algebra and $$(\cdot, ... , \cdot) \colon (\Nil_{R})^{\times(k+1)} \longrightarrow \G_{a, R}$$ is polylinear map. Then, there exists an element $e\in R$ such that 
$$(a_{0}, ... , a_{k}) = e a_{0}...a_{k}.$$
\end{lemma}
\begin{proof}
Consider a morphism $\G_{a, R} \to \Homgr(\Nil_{R}^{\otimes(k+1)}, \G_{a, R})$, that maps an element $e\in \G_{a, R}(A) = A$ to the map $\psi_{e}$, which is given by $\psi_{e}(a_{0},...,a_{k}) = ea_{0}...a_{k}.$ The statement of the lemma is equivalent to the fact that this is an isomorphism. We will prove the statement by induction on $k$. 

If $k = 0$, then a morphism $\Nil_{R} \to \G_{a, R}$ is the same as a power series $f \in R[[x]]$ such that $f(x + y) = f(x) + f(y).$ From this equality, it is clear that $f$ has no constant term. Write $f = ex + a_{k}x^{k} + a_{k+1}x^{k+1}+...$. By comparing the coefficients of $xy^{k-1}$ in the equality $f(x+y) = f(x) + f(y)$, we obtain $a_{k}k = 0$. Since $R$ is an $\Q$-algebra, the base of the induction is proved.

Now, using induction,
$$\Homgr(\Nil_{R}^{\otimes(k+1)}, \G_{a, R}) = \Homgr(\Nil_{R}^{\otimes k}, \Homgr(\Nil_{R}, \G_{a, R})) = \G_{a, R}.$$ 

Moreover, one can see that this identification is induced by the morphism described at the beginning of the proof. The proof is complete.
\end{proof}

\begin{lemma}
\label{nilppolymult}
Let $R$ be a $\Q$-algebra and $$(\cdot, ... , \cdot) \colon (\Nil_{R})^{\times(k+1)} \longrightarrow \G_{m, R}$$
be a polylinear map.
Then
\begin{enumerate}
\item There exists an element $e\in R,$ such that $(a_{0}, ... , a_{k}) = \exp(ea_{0}...a_{k}).$
\item Assume, furthermore, that there exists a collection of integers $s_{0}, ... , s_{k}\in \Z$ such that $\sum s_{i} \neq 0$ and for any $m \in \Z$, the following equality holds: $(a_{0}, ... , a_{k}) = (m^{s_{0}}a_{0}, ... ,m^{s_{k}} a_{k})$. Then, $(\cdot, ..., \cdot) = 1.$
\end{enumerate}
\end{lemma}
\begin{proof}
Let us prove the first part. A morphism to the affine scheme $\G_{m, R}$ is given by an element $f(x_{0}, ... , x_{k})\in \mathcal{O}(\Nil^{\times(k+1)})^{*} = R[[x_{0}, ... , x_{k}]]^{*}$. By polylinearity, it follows that $f(0, ... , 0) = 1.$ From this, we see that the image of $(\cdot, ... , \cdot)$ lies in the subfunctor  $1 + \Nil_{R}.$ Next, using an isomorphism $\log\colon 1 + \Nil_{R} \to \Nil_{R}$, we obtain a polylinear map $\Nil_{R}^{\times(k+1)} \to \Nil_{R}.$ By Lemma \ref{nilppolylinear}, we see that there exists an element $e\in R$ such that $\log((a_{0}, ... , a_{k})) = ea_{0} ... a_{k}.$ Now, using the map $\exp$, we are done.

Next, we prove the second part. Using the first part and the conditions of the lemma, we obtain 
$$\exp(e a_{0} ... a_{k}) = \exp(e m^{\sum s_{i}} a_{0} ... a_{k}).$$

Now, since $R$ is a $\Q$-algebra, we see that $e = 0.$ Lemma is proved.
\end{proof}

Proposition \ref{nilpotents} is a direct consequence of Proposition 5.16 from the work \cite{GorOsi15}. It can be easily derived from the theory of thick ind-cones developed by the authors. Nevertheless, we provide a direct proof of this fact. To begin with, we state a simple lemma.

\begin{lemma}\label{obvious}
Let X be an affine scheme over a ring $R$. Suppose we are given two morphisms $\alpha, \beta \colon \mathbb{A}_{R}^{k} \to X$ that coincide after restriction to $\Nil_{R}^{k}$. Then, $\alpha = \beta.$
\end{lemma}
\begin{proof}
Since $X$ is affine, the morphisms $\alpha, \beta$ are uniquely determined by the induced morphisms $\mathcal{O}(X) \to \mathcal{O}(\mathbb{A}^{k})$. Now, everything follows from the fact that the map $\mathcal{O}(\mathbb{A}_{R}^{k}) \to \mathcal{O}(\Nil_{R}^{k})$ is an inclusion.

\end{proof}

\begin{propos}
\label{nilpotents}
Let $R$ be a ring, and let $\alpha, \beta \colon (\Ladd_{R}^{\#})^{k} \to \G_{m,R}$ be two polylinear morphisms that coincide on $A$-points for all $R$-algebras of the form $(f_{1}, ... , f_{k}),$ where $f_{i} = a_{i}t^{l},$ with $l\in \Z^{n}$ and $a_{i}\in \Nil(A)$. Then, $\alpha = \beta$.
\end{propos}
\begin{proof}
Fixing the arguments, we may assume that $k = 1$. To simplify the notation, we put
$$X_{\lambda} =  \Ladd_{\lambda, R}, \ T_{\lambda} = \O(X_{\lambda}),$$
$$X_{\lambda, N} =  \Ladd_{\lambda, N, R}^{\#}, \ T_{\lambda} = \O(X_{\lambda, N}).$$

Next, for a finite subset $S\subset \Z^{n}$ we denote by $X_{\lambda, N, S}$ the subfunctor of the functor $X_{\lambda, N},$ which assigns to an $R$-algebra $A$ the set of $f\in X_{\lambda, N}(A)$ such that supp$(f) \subset S.$ By Lemma \ref{obvious} and the conditions of the proposition, we obtain that $\alpha = \beta$ after restriction to $X_{\lambda, N, S}.$ Similarly, we denote $T_{\lambda, N,S} = \O(X_{\lambda, N, S}).$

Recall that $T_{\lambda}$ is simply a polynomial algebra generated by variables $x_{l},$ where $l \in \Z^{n}_{\lambda}.$ Now, the algebra $T_{\lambda, N}$ is the quotient of $T_{\lambda}$ by an ideal $I_{\lambda, N}$, which is generated by homogeneous elements of degree not less than $N$.

Using the discussion above, we see that we may equip $T_{\lambda, N}$ with the natural grading, where the degree of each variable $x_{l}$ is $1$. Moreover, the epimorphism $T_{\lambda} \to T_{\lambda, N}$ is an isomorphism in all graded components up to degree $N-1$. Thus, in graded components with degree at most $N-1$, the algebra $T_{\lambda, N}$ is freely generated by monomials of variables $x_{l}, l \in \Z^{n}_{\lambda}.$ Similarly, the algebra $T_{\lambda, N, S}$, in graded components with degree at most $N-1$, is freely generated by monomials in the variables $x_{l}$, where $l\in S.$  

Let $\alpha_{\lambda, N}, \beta_{\lambda, N}$ denote the restrictions of the morphisms $\alpha, \beta$ to $X_{\lambda, N}.$ Then, by the Yoneda lemma, these restrictions are uniquely determined by invertible elements $a_{\lambda, N}$, $b_{\lambda, N}$ respectively. We denote by $(a_{\lambda, N})_{s}, (b_{\lambda, N})_{s}$ the corresponding graded components. It suffices to show that $(a_{\lambda, N})_{s} = (b_{\lambda, N})_{s}$ for ${s \leq N-1}.$ Indeed, for $s \geq N$, we have $(a_{\lambda, s+1})_{s} =  (b_{\lambda, s+1})_{s}.$ Moreover, the elements $a_{\lambda, s+1}, b_{\lambda, s+1}$ map to $a_{\lambda, N}, b_{\lambda, N}$, respectively, under the map $T_{\lambda, s+1} \to T_{\lambda, N}.$

Now, $T_{\lambda, N}$ in degrees up to $N-1$ is generated by monomials in the variables $x_{l}.$ For any such monomial $x_{l_{1}}... x_{l_{r}}$, let $S = \{l_{1}, ... , l_{r}\}.$ Since $\alpha = \beta$ when restricted to $X_{\lambda, N, S},$ the coefficients of the monomial $x_{l_{1}}... x_{l_{r}}$ are the same. From this, the desired result follows.
\end{proof}

\section{Proof of the Main Result} \label{proofs}

\begin{lemma}\label{mainlemma} Let $R$ be a $\Q$-algebra, and let
    
    $$\{\cdot, ..., \cdot\} \colon (\Ladd_{R}^{\#})^{\times(k+1)} \longrightarrow \G_{m, R}$$ 
be a polylinear morphism of functors (defined on the category of $R$-algebras) that is invariant under continuous automorphisms of the form $t_{i} \mapsto m_{i}t_{i}$ where $m_{i} \in \Z$. Then the following holds:
\begin{enumerate}
	\item The image of the morphism lies in the subfunctor $1 + \Nil_{R}.$
	\item Let $A$ be an arbitrary $R$-algebra, and let $a_{0}, ... , a_{k}\in A$, $l_{0}, ... , l_{k} \in \Z^{n}$ be such that $a_{i}t^{l_{i}}\in \Ladd_{R}^{\#} $. Explicitly, this means that if $l_{i} \leq 0,$ then $a_{i}$ is nilpotent. Then there exists an element $e(l_{0}, ... , l_{k}) \in R,$ depending only on $l_{0}, ... , l_{k},$ such that the following formula holds:
$$\log(\{a_{0}t^{l_{0}}, ... , a_{k}t^{l_{k}} \}) = e(l_{0}, ..., l_{k}) a_{0}...a_{k}. $$    
Moreover, if $\sum\limits_{i = 0}^{k} l_{i} \neq 0,$ then $e(l_0, ..., l_k) = 0.$
	\item The map $\{ \cdot, ... , \cdot \}$ is uniquely determined by the constants $e(l_{0}, ... , l_{k})$. Equivalently, if $A$ is an $R$-algebra and $f_{0}, ... , f_{k} \in \Ladd^{\#}(A),$ then the following formula holds:
$$ \{f_{0}, ..., f_{k} \} = \prod\limits_{l_{0} + ... + l_{k} = 0} \exp(e(l_{0}, ..., l_{k}) a_{l_{0}}^{(0)}...a_{l_{k}}^{(k)}),$$
where $f_{i} = \sum a_{l}^{(i)}t^{l}.$
\end{enumerate}
    
\end{lemma}

\begin{proof}
Fix elements $l_{0}, ... , l_{k} \in \Z^{n}.$ Consider the multilinear morphism $\Nil_{R}^{\times(k+1)} \to \G_{a, R},$ that maps $(a_{0}, ... , a_{k})$ to $\{a_{0}t^{l_{0}}, ... , a_{k} t^{l_{k}}\}.$ From Lemma \ref{nilppolymult} it follows that there exists an element $e(l_{0}, ... , l_{k})\in R$ such that the following holds:
\begin{equation} \label{multeq}
\{a_{0}t^{l_{0}}, ... , a_{k} t^{l_{k}}\} = \exp(e(l_{0}, ... , l_{k}) a_{0}...a_{k}).
\end{equation}

From the hypotheses and part 2 of Lemma \ref{nilppolymult}, we conclude that if $\sum l_{i} \neq 0,$ then $e(l_{0}, ... , l_{k}) = 0.$

We now consider the morphism of functors $\varphi \colon \Ladd_{R}^{\times(k+1)} \to \G_{m, R}$, defined by the formula:
\begin{equation}\label{bigprod}
\varphi(f_{0}, ..., f_{k}) = \prod\limits_{l_{0} + ... + l_{k} = 0} \exp(e(l_{0}, ..., l_{k}) a_{l_{0}}^{(0)}...a_{l_{k}}^{(k)}),
\end{equation}
where $f_{i} = \sum a_{l}^{(i)}t^{l}.$ Note that this is well-defined because for any collection $f_{0}, ... , f_{k} \in \Ladd^{\#}(A)$, there exists $\lambda \in \Lambda_{n},$ such that supp$(f_{i})\subset \Z^{n}_{\lambda}$, and the fibers of the summation map
$$\Z^{n}_{\lambda} \times ... \times \Z^{n}_{\lambda} \longrightarrow \Z^{n}$$
are finite. Therefore, the product in \eqref{bigprod} is finite for any fixed collection $f_{0}, ... , f_{k} \in \Ladd^{\#}(A)$.

From the arguments above, it follows immediately that $\varphi$ and $\{\cdot, ... , \cdot \}$ coincide on collections $f_{0}, ... , f_{k} \in \Ladd^{\#}(A),$ where each $f_{i}$ has the form $at^{l}, a \in \Nil(A).$ Then by Proposition \ref{nilpotents} we conclude that they coincide identically. This proves the lemma.
\end{proof}

Let $S \subset \{1, ... , n\}$ and $K\subset \{1, ... , n\}$ be a complementary subsets, i.e., $S \sqcup K = \{1, ... , n\}$. If $\Z^{n}$ the free abelian group with basis $e_{i},$ then $S$ and $K$ induce a decomposition $\Z^{n} = \Z^{S} \oplus \Z^{K}.$ Explicitly, $\Z^{S}$ is generated by $e_{i}, i \in S$, and $\Z^{K}$ is generated by $e_{i}, i \in K$. Let $P_{S}, P_{K}\colon \Z^{n} \to \Z^{n}$ denote the projections onto $\Z^{S}$ and $\Z^{K}$, respectively. Furthermore, let us denote $|S|$ and $|K|$ the cardinalities of these sets. Then we identify $\Z^{S}$ with $\Z^{|S|}$ using the naturally ordered basis $e_{i}, i \in S$. Analogously, we identify $\Z^{K}$ with $\Z^{|K|}$.

\begin{propos}
\label{mainprop}
Let $R$ be a $\Q$-algebra, and let    
$$\{\cdot, ..., \cdot\} \colon (\Ladd_{R}^{\#})^{\times(k+1)} \longrightarrow \G_{m, R}$$
be a polylinear morphism of functors (defined on the category of $R$-algebras) that is invariant under continuous automorphisms of the form $t_{i} \mapsto m_{i}t_{i}$ where $m_{i} \in \Z$. Assume further that $S \sqcup K = \{1, ... , n\}$ are sets such that $|S| = k \leq n $ and $\{\cdot, ... , \cdot\}$ is invariant under continuous automorphisms $\varphi \in Aut(\mathcal{L}^{n}(\Q))$ satisfying $\varphi(t_{i}) =  t_{i}$ for $i \in K$. Then

\begin{enumerate}
\item The image of the morphism lies in the subfunctor $1 + \Nil_R.$
\item There exists $e\in R,$ such that for any $R$-algebra $A$ and arbitrary elements $f_{0}, ... , f_{k} \in \Ladd_{R}^{\#}$, the following holds: 
$$ \{f_{0}, ..., f_{k} \} = \prod\limits_{l_{0} + ... + l_{k} = 0} \exp(e \det(v_{1}, ... , v_{k}) a_{l_{0}}^{(0)}...a_{l_{k}}^{(k)}),$$
where $f_{i} = \sum a_{l}^{(i)}t^{l},$ and $v_{i} = P_{S}(l_{i}) \in \Z^{k}.$
\item If there exist other subsets $S' \neq S$ and $K' \neq K$ satisfying the conditions of the proposition, then $\{\cdot, ... , \cdot\} = 1.$
\end{enumerate}
\end{propos}
\begin{proof}
The first part is contained in Lemma \ref{mainlemma}. Let us prove the second part. By Lemma \ref{mainlemma}, there exist elements $e(l_{0}, ... , l_{k}) \in R$ such that 
\begin{equation}
\label{longprod}
\{f_{0}, ..., f_{k} \} = \prod\limits_{l_{0} + ... + l_{k} = 0} \exp(e(l_{0}, ..., l_{k}) a_{l_{0}}^{(0)}...a_{l_{k}}^{(k)}),
\end{equation}
where $f_{i} = \sum a_{l}^{(i)}t^{l}.$ 

Obviously, it suffices to show that there exists $e\in R$ such that for any $l_{0}, ... , l_{k}\in \Z^{n}$ satisfying $l_{0} + ... + l_{k}= 0$, the following holds: 
$$e(l_{0}, ... , l_{k}) = e\det(v_{1}, ... , v_{k}),$$
where $v_{i} = P_{S}(l_{i}) \in \Z^{k}$. Put also $w_{i} = P_{K}(l_{i})$. Let $TR = R[\eta],$ where $\eta^{k+2} = 0.$ Then by formula \eqref{longprod}, for any $f_{0}, ... , f_{k} \in \Ladd(R)$, the following holds:
\begin{equation}
\label{eqeta}
\log \{\eta f_{0}, ...,\eta f_{k} \} = \sum\limits_{l_{0} + ... + l_{k} = 0} e(l_{0}, ..., l_{k}) a_{l_{0}}^{(0)}...a_{l_{k}}^{(k)}\eta^{k+1}.
\end{equation}

Note that the right-hand side of the formula lies in $\eta^{k+1}R \subset TR.$ This allows us to define a map:
$$(\cdot, ... , \cdot) \colon (\Ladd(R))^{\times(k+1)} \longrightarrow \G_{a}(R)$$
using the following relation: 
$$\log(\{\eta f_{0} t^{w_{0}} , ... , \eta f_{k} t^{w_{k}}\}) = (f_{0}, ... , f_{k})\eta^{k+1}.$$

It is immediately verified that the map $(\cdot, ... , \cdot)$ is invariant under continuous automorphisms from $Aut(\mathcal{L}^{n}(\Q))$. The continuity and $R$-polylinearity follow from formula \eqref{eqeta}. Then by Theorem \ref{theorem1}, there exists $e'(w_{0}, ... , w_{k})$ such that:
\begin{equation}
\label{eqres}
(f_{0}, ... , f_{k}) = e'(w_{0}, ... , w_{k}) \Res(f_{0} df_{1} \wedge ... \wedge df_{k}).
\end{equation}

Now, combining equations \eqref{eqeta} and \eqref{eqres}, we obtain:
$$e'(w_{0}, ... , w_{k}) \det(v_{1}, ... , v_{k})= e(l_{0}, ... , l_{k}).$$

We need to prove that $e'(w_{0}, ... , w_{k})$ is independent of the $w_{i}$. Let $S = \{i_{1}, ... ,i_{k}\},$ where $i_{1} < i_{2} < ... < i_{k}.$ Set $y_{s} = t_{i_{s}}$ for $1 \leq s \leq k,$ and $y_{0} = (y_{1}...y_{k})^{-1}$. It is clear that it suffices to prove the equality:
$$\{ \eta y_{0}, ... , \eta y_{k}\} = \{ \eta y_{0} t^{w_{0}}, ... , \eta y_{k} t^{w_{k}} \}.$$ 

This follows by repeated application of Lemma \ref{maintech}. Let us now prove the third part. Let $e'\in R$ be the constant from part 2 for the new subsets $S', K'.$ For any collection $l_{0}, ... , l_{k}$ satisfying $l_{0} + ... +l_{k} = 0$, denote $v_{i} = P_{S}(l_{i}), v_{i}' = P_{S'}(v_{i}).$ Then by what has been proved $$e(l_{0}, ... , l_{k}) = e \det(v_{1}, ... , v_{k}) = e(v_{0}, ... , v_{k}) = e' \det(v_{1}', ... , v_{k}') = 0,$$ since $\det(v_{1}', ... , v_{k}')$ has a zero row.
\end{proof}

\begin{cons}
Let $R$ be a $\Q$-algebra, and let    
$$(\cdot, ..., \cdot) \colon (\Lmult_{R}^{\#})^{\times(n+1)} \longrightarrow \G_{m, R}$$
be a polymultiplicative morphism of functors that is invariant under continuous automorphisms from $Aut(\mathcal{L}^{n}(\Q))$. Then there exists an element $e\in R$ such that for any $R$-algebra $A$ and elements $f_{0}, ... , f_{n} \in \Lmult^{\#}(A)$, the following formula holds:
$$(f_{0}, ... , f_{n}) = \exp(e \Res( \log f_{0} \frac{df_{1}}{f_{1}} \wedge ... \wedge  \frac{df_{n}}{f_{n}})).$$
\end{cons}
\begin{proof}
Consider the morphism
$$\{\cdot, ..., \cdot\} \colon (\Ladd_{R}^{\#})^{\times(n+1)} \longrightarrow \G_{m, R}$$
defined by the formula:
$$\{f_{0}, ... , f_{n} \} = (\exp(f_{0}), ... , \exp(f_{n})).$$

It is easily seen that $\{ \cdot, ... , \cdot \}$ satisfies Proposition \ref{mainprop} for $k = n$ and $K = \emptyset.$ Then by part 2, there exists an element $e$ such that: 

$$ \{f_{0}, ..., f_{n} \} = \prod\limits_{l_{0} + ... + l_{n} = 0} \exp(e \det(l_{1}, ... , l_{n}) a_{l_{0}}^{(0)}...a_{l_{n}}^{(n)}) = \exp(e\Res(f_{0} df_{1}\wedge ... \wedge df_{n})),$$
where $f_{i} = \sum a_{l}^{(i)}t^{l}.$ 
Now observe that $(f_{0}, ... , f_{n}) = \{ \log(f_{0}), ... , \log(f_{n}) \}. $Therefore,
$$(f_{0}, ... , f_{n}) = \exp(e \Res( \log f_{0} d\log (f_{1}) \wedge ... \wedge  d\log(f_{n}))).$$

Furthermore, it is straightforward to verify that $d\log(f) = \frac{df}{f}.$ The proof is complete.
\end{proof}

\begin{lemma}
\label{Lmulttrivial}
Let $R$ be a ring, and let 
$$(\cdot, ..., \cdot) \colon \Lmult_{R}^{\times(n+1)} \longrightarrow \G_{m, R}$$
be a polymultiplicative morphism of functors. If $A$ is an $R$-algebra and ${a\in A^{*}}$, $ f_{2}, ..., f_{n}\in\Lmult(A),$ then
\begin{enumerate}
\item For any $b\in A^{*}$ the equality $(a, b, f_{2}, ..., f_{n}) = 1$ holds. 
\item For any $g \in (\Lmult)^{\#}(A)$ the equality $(a, g, f_{2}, ..., f_{n}) = 1$ holds.
\end{enumerate}
\end{lemma}
\begin{proof}
We prove the first part. Fix $f_{2}, ... , f_{n}$ and consider the induced morphism
$$\G_{m, A} \times \G_{m,A} \longrightarrow \G_{m, A}$$
that maps a pair $(a,b)$ to $(a, b, f_{2}, ... , f_{n}).$This morphism is bimultiplicative and therefore trivial by part 2 of Corollary \ref{trivial}. The second part now follows analogously using part 4 of Corollary \ref{trivial}.
\end{proof}

\begin{propos}\label{propdet}
Let $R$ be a $\Q$-algebra, and let
$$(\cdot, ... \cdot) \colon \Lmult^{\times(n+1)}_{R} \longrightarrow \G_{m, R} $$
be a polymultiplicative morphism of functors that is invariant under continuous automorphisms from $Aut(\mathcal{L}^{n}(\Q))$. Then there exists an element $m \in \underline{\Z}(R),$ such that for any $R$-algebra $A$ and elements $f_{1}, ... , f_{n} \in \Lmult(A)$, the following holds:
\begin{equation}
\label{det}
(a, f_{1}, ... , f_{n}) = a^{m\det(v(f_{1}), ... , v(f_{n}))}.
\end{equation}
Moreover, for every $1 \leq i \leq n$, the following relation holds:
\begin{equation}
\label{antisymmetry}
(a, f_{1}, ..., f_{i} , ... , f_{n}) = (f_{i}, f_{1}, ..., a , ... , f_{n})^{-1}. 
\end{equation}
\end{propos}

\begin{proof}
Let $1 \leq i_{1}, ... , i_{n} \leq n $ be an arbitrary set of indices. We first prove that if this collection contains repeated indices, then for every $a\in \NN \subset R$, the following identity holds:
\begin{equation}
\label{rel1}
(a, t_{i_{1}}, ... , t_{i_{n}}) = 1.
\end{equation}

Since there are two repeated indices, by the pigeonhole principle, there exists $1 \leq j \leq n$ such that $j$ does not appear in the set $\{i_{1}, ... , i_{n}\}$. There exists an automorphism that maps $t_{j}$ to  $at_{j}$ and $t_{i}$ to $t_{i}$ for  $i \neq j.$ Using the invariance condition, we obtain:
$$(t_{j}, t_{i_{1}}, ... , t_{i_{n}}) = (at_{j}, t_{i_{1}}, ... , t_{i_{n}}).$$

Therefore, we conclude that $(a, t_{i_{1}}, ... , t_{i_{n}}) = 1.$

Now, let $1 \leq i_{1}, ... , i_{n} \leq n $ be a set of distinct indices. We now show that for all integers $1 \leq k, s \leq n$ and $a\in \NN$ the following holds:
\begin{equation}
\label{rel2}
 (a, t_{i_1}, ...,  t_{i_s}, ..., t_{i_k}, ..., t_{i_n}) = (a, t_{i_1}, ...,  t_{i_k}, ..., t_{i_s}, ..., t_{i_n})^{-1}.
\end{equation}

Without loss of generality, assume $i_{k} > i_{s}$. There exists an automorphism that maps $t_{i_{k}}$ to $t_{i_{k}}t_{i_{s}}$ and $t_{i_{l}}$ to $t_{i_{l}}$ for $l \neq k$. By the invariance condition, we obtain:
$$ (a, t_{i_1}, ...,  t_{i_k}, ..., t_{i_k}, ..., t_{i_n}) = (a, t_{i_1}, ...,  t_{i_{k}}t_{i_{s}}, ..., t_{i_{k}}t_{i_{s}}, ..., t_{i_n}).$$

Expanding the brackets and applying identity \eqref{rel1}, we obtain the required result.

We now prove the proposition. We define a morphism $\G_{m,R} \to \G_{m,R},$ which maps an element $a$ to the value $(a, t_{1}, ... , t_{n})$. By Lemma \ref{Gm} there exists an element $m\in \underline{\Z}(R),$ such that $(a, t_{1}, ... , t_{n}) = a^{m}.$ 

Now, let $f_{1}, ... , f_{n}$ be elements of $\Lmult(A)$. As explained in Section \ref{maindef}, any $f_{i}$ admits the decomposition $f_{i} = a_{i}g_{i} t^{v(f_{i})},$ where $a_{i} \in A^{*}$, ${g_{i} \in \Lmult^{\#}(A)}$, $v(f_{i}) \in \underline{\Z}^{n}(A).$

Applying Lemma \ref{Lmulttrivial} and the polymultiplicativity, we derive:
$$(a, f_{1}, ... , f_{n}) = (a, t^{v(f_{1})}, ... , t^{v(f_{n})}).$$

It is straightforward to show that any polylinear morphism of functors ${(\underline{\Z}_{R}^{n})^{\times n} \to \underline{\Z}_{R}}$ is automatically $\underline{\Z}_{R}$-polylinear. Combining relations \eqref{rel1} and \eqref{rel2} we obtain that for any $a\in \NN$:
$$(a, f_{1}, ... , f_{n}) = (a, t_{1}, ... ,t_{n})^{\det(v(f_{1}), ... , v(f_{n}))} = a^{m\det(v(f_{1}), ... , v(f_{n}))}.$$

By functoriality, this identity holds for all ${a\in A^{*}}.$ It remains to verify the antisymmetry relation \eqref{antisymmetry}. It is easily seen that it suffices to consider the case when $f_{1} = t_{1}, ... , f_{n} = t_{n}, a\in \NN$. Consider the automorphism that maps $t_{i}$ to  $at_{i}$ and $t_{j}$ to  $t_{j}$ for $j \neq i.$ We then obtain the equality: 
$$(t_{i}, t_{1}, ..., t_{i}, ... , t_{n}) = (at_{i}, t_{1}, ... ,at_{i}, ... , t_{n}). $$

Taking into account relation \eqref{rel1} and Lemma \ref{Lmulttrivial}, we obtain the desired result.
\end{proof}

\begin{propos}\label{sign}
Let $R$ be a ring. 
\begin{enumerate}
\item For any map $\omega \colon (\underline{\Z}^{n}_{R})^{\times(n+1)} \longrightarrow \G_{m,R},$
that is invariant under the action of the group $B(n, \underline{\Z}({R})),$ one can define
$$(\cdot, ... , \cdot)_{\omega} \colon \Lmult^{\times(n+1)}_{R} \longrightarrow \G_{m, R}$$
via the formula $(f_{0}, ... , f_{n})_{\omega} = \omega(v(f_{0}), ... , v(f_{n}))$. Then $(\cdot, ... , \cdot)_{\omega}$ is polymultiplicative and invariant under continuous automorphisms.
\item Let $$(\cdot, ... , \cdot) \colon \Lmult^{\times(n+1)}_{R} \longrightarrow \G_{m, R}$$
be a polymultiplicative morphism of functors that is invariant under continuous automorphisms from $Aut(\mathcal{L}^{n}(\Q))$. Then one can define $\omega \colon (\underline{\Z}^{n}_{R})^{\times(n+1)} \to \G_{m,R}$ by:
$$\omega(l_{0}, ... , l_{n}) = (t^{l_{0}}, ... , t^{l_{n}}). $$
This morphism will be invariant under the action of $B(n, \underline{\Z}({R}))$and polylinear. 
\end{enumerate}
\end{propos}
\begin{proof}
The first claim follows directly from equality \eqref{eqauto}. Indeed,
$$(\varphi(f_{0}), ... , \varphi(f_{n}))_{\omega} = \omega(M(\varphi)(v(f_{0})), ..., M(\varphi)(v(f_{n}))) = \omega(v(f_{0}), ... , v(f_{n})).$$

Let us prove the second part. Let $U \in B(n, B(n, \underline{\Z}({R}))).$ Then
$$\omega(U(l_{0}), ... , U(l_{n})) = (t^{U(l_{0})}, ... , t^{U(l_{n})}) = (U(t^{l_{0}}), ... , U(t^{l_{n}}) = (t^{l_{0}}, ... , t^{l_{n}}) . $$
\end{proof}
\begin{remark}
If $R= \Z$, then the image of such a morphism always lies in $\{\pm 1\}$.
\end{remark}
\begin{remark}
Let $n = 2$, and let $e_{1}, e_{2}$ be a basis of the free abelian group $\Z^{2}$. Denote $a_{ijk} = \omega(e_{i}, e_{j}, e_{k}) \in R^{*}$ for $1 \leq i, j, k \leq 2.$ Solving the corresponding system of equations yields the following set of conditions: $a_{111} = a_{112}^{2} = 1$, $a_{112} = a_{211} = a_{121}$, $a_{221}a_{212}a_{122} = a_{112}$ and $a_{222}\in R^{*}$ can be any invertible element.
\end{remark}

\begin{propos}\label{finalprop}
Let $R$ be a $\Q$-algebra, and let
$$(\cdot, ... \cdot) \colon \Lmult^{\times(n+1)}_{R} \longrightarrow \G_{m, R} $$
be a polymultiplicative morphism of functors that is invariant under continuous automorphisms from $Aut(\mathcal{L}^{n}(\Q))$.  By Proposition \ref{propdet} there exists an element $m\in \underline{\Z}(R)$ such that formula \eqref{det} holds. Then, for any $R$-algebra $A$ and for any $f_{0} \in \Lmult^{\#}(A)$, $f_{1}, ... , f_{n} \in  \Lmult(A)$, the following equality holds:

$$ (f_{0}, ... , f_{n}) = CC_{n}(f_{0}, ... , f_{n})^{m}.$$
\end{propos}
\begin{proof}
Let $1 \leq i_{1}, ... , i_{s} \leq n$.  We prove that there exists an element $e\in R$ such that for any $f_{0}, ... , f_{k} \in \Lmult^{\#}(A)$, where $k + s = n,$ the following holds:
\begin{align}\label{bigformula}
(f_{0}, ... , f_{k}, t_{i_{1}}, ... , t_{i_{s}}) = \exp(e \Res(\log(f_{0}) \omega_{1} \wedge \omega_{2})),
\end{align}
where $\omega_{1} = \frac{df_{1}}{f_{1}} \wedge ... \wedge \frac{df_{k}}{f_{k}}$, $\omega_{2} = \frac{dt_{i_{1}}}{t_{i_{1}}} \wedge ... \wedge \frac{dt_{i_{s}}}{t_{i_{s}}}.$

First, we show how the proposition follows from formula \eqref{bigformula}. Let $m\in \underline{\Z}(R)$ be as in Proposition \ref{propdet} and let $f_{1}, ... , f_{n} \in \Lmult(A)$. Each element $f_{i}$ can be decomposed as a product $f_{i} = a_{i}g_{i}t^{v(f_{i})},$ where $a_{i} \in \G_{m}(A)$, $g_{i} \in \Lmult^{\#}(A).$ By using multimultiplicativity and functoriality, it is easy to see that it suffices to prove the proposition for tuples $f_{0}, ... , f_{n}$ where each $f_{i}$ is either in $\G_{m}(A),$ in $\Lmult^{\#}(A)$, or of the form $t_{j}$ for some $j$. If at least one $f_{i}$ belongs to $\G_{m}(A),$  claim follows directly from Proposition \ref{propdet} and equality \eqref{antisymmetry}. 

Now suppose that any $f_{i}$ is either belongs to $\Lmult^{\#}(A)$ or $f_{i} = t_{j}.$ Without loss of generality, we may assume $f_{1}, ... , f_{k} \in \Lmult^{\#}(A)$, ${f_{k+1} = t_{i_{1}}}, ... , {f_{n} = t_{i_{s}}}$, $k + s = n.$ If among $i_{1}, ... , i_{s}$ there are two repeated indices, then the proposition follows from the above formula and the explicit formula for $CC_{n}.$ 

If all $i_{0}, ... , i_{s}$ are distinct, let $j_{1}, ... , j_{k}$ be the complement of $\{i_{1}, ... , i_{s} \}$ in $\{1, ... , n \}, \ k + s = n$. Then substituting into the formula: $f_{0} \in 1 + \Nil(A), f_{1} = \exp(a_{1} t_{j_{1}}), ... ,\\ f_{k} = \exp(a_{k} t_{j_{k}}),$ where $a_{i} \in \Nil(A),$ we obtain 
$$(f_{0}, ... , f_{k}, t_{i_{1}}, ... , t_{i_{k}}) = \exp(e\varepsilon f_{0}a_{1}...a_{k}),$$
where $\varepsilon \in \{\pm \}$ is determined by the condition $ dt_{j_{1}} \wedge ... \wedge dt_{j_{k}} \wedge dt_{i_{1}} \wedge ... \wedge dt_{i_{s}} = \varepsilon dt_{1} \wedge ... \wedge dt_{n}.$

Now from Proposition \ref{propdet} it follows that
$$(f_{0}, ... , f_{k}, t_{i_{1}}, ... , t_{i_{s}}) = \exp(m\varepsilon f_{0}a_{1}...a_{k}).$$

This shows that $e = m,$ and the proposition follows from the explicit formula for $CC_{n}.$

We now prove the formula \eqref{bigformula}. Define
$$\{ \cdot, ... , \cdot \} \colon ((\Ladd)^{\#}_{R})^{\times(k+1)} \longrightarrow \G_{m,R} $$
by a formula
$$\{ g_{0}, ... , g_{k} \} = (\exp(g_{0}), ... ,\exp(g_{k}), t_{i_{1}}, ... , t_{i_{s}}).$$

Let $S \subset \{ 1, ... , n\}$ be a subset with $|S| = k$ such that $S \cap \{ i_{1}, ... , i_{s}\} = \emptyset,$ and let $K$ be the complement of $S$ in $\{1, ... , n\}.$ Note that we always have the inclusion $\{i_{1}, ... , i_{s}\} \subset K$ with an equality precisely when all $i_{1}, ... , i_{s}$ are distinct. Then by part 2 of Lemma \ref{Lmulttrivial}, the morphism $\{ \cdot, ... ,\cdot \}$ satisfies Proposition \ref{mainprop} with chosen $S, K$. Moreover, if the turple $i_{1}, ... , i_{s}$ contains repeated elements, one can choose another such pair $S, K$. In this case, part 3 of Proposition \ref{mainprop} immediately yields: $(f_{0}, ... , f_{k}, t_{i_{1}}, ... , t_{i_{s}}) = 1.$

Now suppose the elements $i_{1}, ... , i_{s}$ are all distinct. Then the following holds: 
$$\Res(t^{l_{0}} dt^{l_{1}} \wedge ... \wedge dt^{l_{k}} \wedge \frac{dt_{i_{1}}}{t_{i_{1}}} \wedge ... \wedge \frac{dt_{i_{s}}}{t_{i_{s}}}) = \pm \det(v_{1}, ... , v_{k}),$$
where the elements $v_{i}$ were defined previously. This immediately implies that
$$\{\log(f_{0}), ... , \log(f_{k})\} = \exp(\pm e \Res(\log(f_{0}) \frac{df_{1}}{f_{1}} \wedge ... \wedge \frac{df_{k}}{f_{k}} \wedge \frac{dt_{i_{1}}}{t_{i_{1}}} \wedge ... \wedge \frac{dt_{i_{s}}}{t_{i_{s}}})).$$

Now, using the definition of $\{\cdot, ... , \cdot\}$ and replacing $e$ with -$e$ if necessary, we obtain the desired result.
\end{proof}

\begin{theorem}
\label{mainthm}
Let $R$ be a $\Q$-algebra, and let
$$(\cdot, ... \cdot) \colon \Lmult^{\times(n+1)}_{R} \longrightarrow \G_{m, R} $$
be a polymultiplicative morphism of functors that is invariant under continuous automorphisms from $Aut(\mathcal{L}^{n}(\Q))$. By Proposition \ref{propdet}, there exists an element $m\in \underline{\Z}(R)$ for which formula \eqref{det} holds. Then, there exists a morphism $\omega \colon (\underline{\Z}^{n}_{R})^{\times(n+1)} \to \G_{m,R}$ that is invariant under the action of the group $B(n, \underline{\Z}(R))$, and satisfies the following condition:
\begin{equation}\label{eqwithomega}
(f_{0}, ... , f_{n}) = \omega(v(f_{0}), ... , v(f_{n})) CC_{n}(f_{0}, ... , f_{n})^{m}.
\end{equation}
\end{theorem}
\begin{proof}
Define $\omega(l_{0}, ... , l_{n}) = (t^{l_{0}}, ... , t^{l_{n}}) CC_{n}(f_{0}, ... , f_{n})^{-m}$. We now prove the required equality holds.

If $f_{0} \in \G_{m}(A)$, this reduces exactly to Proposition \ref{propdet}. Analogously, if any $f_{i} \in \G_{m}(A),$ we apply equality \eqref{antisymmetry} and Proposition $\ref{propdet}$. Next, when all $f_{i}$ are of the form $t^{l_{i}},$the statement follows directly from the definition of $\omega.$

If $f_{0} \in \Lmult^{\#}(A),$ the statement follows directly from Proposition \ref{finalprop}. When for some $i\neq 0$ the element $f_{i}$ belongs to $\Lmult^{\#}(A)$ we introduce a new map $$(\cdot, ... , \cdot)'\colon (\Lmult)^{\times(n+1)}_{R} \longrightarrow \G_{m, R},$$
defined by the rule $$(f_{0}, ... f_{n})' = (f_{i}, f_{1}, ...,f_{i-1}, f_{0}, f_{i+1}, ... , f_{n}).$$ 

Therefore, by applying Proposition \ref{finalprop} and the equality \eqref{antisymmetry}, to the map $(\cdot, ..., \cdot)',$ we immediately conclude that $(f_{0}, ... ,f_{n}) = CC_{n}(f_{0}, ... , f_{n})^{m}.$ 

Since every element $f \in \Lmult(A)$ decomposes as $f = a g t^{v(f)},$ where $a \in \G_{m}(A)$ and $g \in \Lmult^{\#}(A),$ the general case reduces to those already considered. This completes the proof of the theorem.
\end{proof}
\begin{remark} The proof uses only automorphisms mapping each $t_{i}$ to a rational function in $t_{1}, ... , t_{n}.$
\end{remark}

\begin{cons} \label{cons}
If additionally the equalities $\{a, t_{1}, ... , t_{n}\} = a$ and $(t^{l_{0}}, t^{l_{1}}, t^{l_{2}}, ... , t^{l_{n}}) = (-1)^{\sgn(l_{0}, ... , l_{n})}$ are satisfied, then the map coincides with $CC_{n}$, i.e. $(\cdot, ... ,\cdot) = CC_{n}.$
\end{cons}
\begin{cons}
Let $R$ be an arbitrary torsion-free ring and let
$$ (\cdot, ... , \cdot)\colon (\Lmult)^{\times(n+1)}_{R} \longrightarrow \G_{m, R}$$
be a polymultiplicative morphism of functors that, when restricted to $R\otimes_{\Z} \Q$-algebras,  satisfies the conditions of Theorem \ref{mainthm}. Then there exist an element ${m \in \underline{\Z}(R)}$ and a morphism $\omega\colon (\underline{\Z}^{n}_{R})^{\times(n+1)} \to \G_{m,R}$ which is invariant under the action of $B(n, \underline{\Z}(R))$, such that the following equality holds:
$$(f_{0}, ... , f_{n}) = \omega(v(f_{0}), ... , v(f_{n})) CC_{n}(f_{0}, ... , f_{n})^{m}.$$
\end{cons}
\begin{proof}
Let $R_{\Q} = R \otimes_{\Z} \Q,$ and denote by $\{\cdot, ..., \cdot\}_{\Q}$ the restriction of $\{\cdot, ... , \cdot\}$ to $R_{\Q}$-algebras. Applying Theorem \ref{mainthm} to $\{\cdot, ... , \cdot \}_{\Q}$, we obtain an element $m \in \underline{\Z}(R_{\Q})$ and a morphism $\omega\colon (\underline{\Z}^{n}_{R_{\Q}})^{\times(n+1)} \to \G_{m,R_{\Q}},$ satisfying equation \eqref{eqwithomega}.

Define a morphism $\G_{m, R} \to \G_{m,R}$ by $a \mapsto \{a, t_{1}, ... , t_{n}\}.$ By Proposition \ref{Gm}, there exists $m'\in \Z(R)$ such that $\{a, t_{1}, ... , t_{n}\} = a^{m'}.$ For $a \in \Z,$ we have $a^{m'} = a^{m}$ in $R_{\Q}.$ Since the canonical map $R\to R_{\Q}$ is injective, $m'$ maps to $m$ under $\underline{\Z}(R) \to \underline{\Z}(R_{\Q}).$ Note that the morphism $\omega$ extends to $R$-algebras if and only if its values on basis elements lie in $R^{*} \subset (R_{\Q})^{*}.$ This condition holds here, so $\omega$ admits an extension. By the uniqueness of extensions from $R_{\Q}$-algebras to $R$-algebras (see \cite[Theorem 6.12]{GorOsi15}) equation \eqref{eqwithomega} holds for all $R$-algebras. This completes the proof.
\end{proof}

\begin{remark}
In the work \cite{GorOsi15}, the multidimensional Contou-Carrère symbol was defined by means of successive application of the boundary map in Quillen $K$-theory: $\delta_{n+1}\colon K_{n+1}(\mathcal{L}(A)) \to K_{n}(A).$ Specifically, the symbol is constructed as the composition 
$$(\cdot, ... , \cdot) \colon \Lmult^{\times(n+1)} \to \mathcal{L}^{n}(K_{n+1}^{M}) \to \mathcal{L}^{n}(K_{n+1}) \to ... \to K_{1} \to \G_{m}.$$

Here, $K_{n}^{M}$ denotes the Milnor $K$-group, i.e., the functor that assigns to a ring $A$ the quotient group $(A^{*})^{\otimes n}$ modulo the subgroup generated by elements $a_{1} \otimes ... \otimes a_{n}$ where $a_{i} + a_{i+1} = 1.$ for some $i$. The first map is the natural projection, the second arises from the Loday product in Quillen $K$-theory,and all subsequent maps are boundary homomorphisms. A highly nontrivial problem is to find an explicit formula for this composite map. The authors address this by computing the tangent space to Milnor $K$-groups (see \cite{GorOsiTangent}). From Corollary \ref{cons} of Theorem \ref{mainthm}, one can derive an alternative proof of the explicit formula.

First, the explicit construction of $\delta_{n+1}$ implies its invariance under $A$-linear automorphisms $\psi\colon A((t)) \to A((t))$, satisfying ${\psi(t) = \sum\limits_{l\geq 1} a_{l}t^{l}}$ (see \cite[Lemma 7.3]{OsipovZhu2}). That is, for such $\psi$ with induced automorphism $\psi_{*}$ on $K_{n+1}(A((t)))$, we have $\delta_{n+1} \circ \psi_{*} = \delta_{n+1}.$ Now let $\varphi \in Aut(\mathcal{L}^{n}(\Q)).$ It admits a decomposition $\varphi = \varphi_{n} \circ ... \circ \varphi_{1},$ where $\varphi_{i}(t_{j}) = t_{j}$ for $i\neq j.$ From this, one immediately deduces the invariance of $(\cdot, ... , \cdot)$ under $Aut(\mathcal{L}^{n}(\Q)).$

The equality $(a, t_{1}, ... , t_{n}) = a$ follows directly from properties of $\delta_{m+1}$ (see \cite[Remark 8.9]{GorOsi15}). Meanwhile, the identity $(t^{l_{0}}, ... , t^{l_{n}}) = (-1)^{\sgn(l_{0}, ... , l_{n})}$ is a consequence of the weak stability property of the ring $\mathcal{L}^{n}(A)$ (see \cite[Step 3 in the proof of Theorem 8.17]{GorOsi15}).
\end{remark}

\vspace{0.3cm}

\noindent National Research University Higher School of Economics, Laboratory of Mirror Symmetry,  6 Usacheva str., Moscow 119048, Russia

\noindent {\it E-mail:}  $levashev.va@phystech.edu$


\begin{thebibliography}{99}
\bibitem{AndersonRomo}
G. Anderson, F. Pablos Romo, {\em Simple proofs of classical explicit reciprocity laws on curves using determinant groupoids over an {A}rtinian local ring,} Comm. Algebra, {\bf 32} (2004), 79--102.

\bibitem{ArbConKac}
E. Arbarello, C. De Concini, V. G. Kac, {\em The Infinite Wedge Representation and the Reciprocity Law for Algebraic Curves,} Theta functions–Bowdoin 1987, Part 1, Proc. Sympos. Pure Math., {\bf 49}, Part 1, Amer. Math. Soc., Providence, RI, 1989, 171--190.

\bibitem{BBE}
A. Beilinson, S. Bloch, H. Esnault, {\em {$\epsilon$}-factors for {G}auss-{M}anin determinants,} Mosc. Math. J., {\bf 2} (2002), 477--532.

\bibitem{CC_{1}}
C. Contou-Carr{\`e}re, {\em Local jacobian, universal {W}itt bivector group and the tame symbol,} C. R. Acad. Sci., Paris, S{\'e}r. I, {\bf 318} (1994), 743--746.

\bibitem{CC2}
C. Contou-Carr{\`e}re, {\em Jacobienne locale d'une courbe formelle relative,} Rendiconti del Seminario Matematico della Universit\`a di Padova, {\bf 130} (2013), 1--106.

\bibitem{Deligne}
P. Deligne, {\em Le symbole mod\'er\'e,} Publ. Math. IHES, {\bf 73} (1991), 147--181.

\bibitem{GorOsi15}
S.O. Gorchinskiy, D.V. Osipov, {\em A~higher-dimensional Contou-Carr\`ere symbol: local theory,} Sb. Math., {\bf 206} (2015), 1191--1259.

\bibitem{GorOsiTangent}
S.O. Gorchinskiy, D.V. Osipov, {\em Tangent space to Milnor K-groups of rings,} Proc. Steklov Inst. Math., {\bf 290} (2015), 34--42.

\bibitem{GorOsi20short}
S.O. Gorchinskiy, D.V. Osipov, {\em Explicit formula for the higher-dimensional Contou-Carr\`ere symbol,} Russian Math. Surveys, {\bf 70} (2015), 171--173.

\bibitem{GorOsi16}
S.O. Gorchinskiy, D.V. Osipov, {\em Higher-dimensional Contou-Carrère symbol and continuous automorphisms,} Funct. Anal. Appl., {\bf 50} (2016), 268--280.

\bibitem{GorOsiCont}
S.O. Gorchinskiy, D.V. Osipov, {\em Continuous homomorphisms between algebras of iterated Laurent series over a ring,} Proc. Steklov Inst. Math., {\bf 294} (2016), 47--66.

\bibitem{GorOsi20}
S.O. Gorchinskiy, D.V. Osipov, {\em Iterated Laurent series over rings and the Contou-Carr\`ere symbol,} Russian Math. Surveys, {\bf 75} (2020), 995--1066.

\bibitem{Kho06}
A. G. Khovanskii, {\em An analog of determinant related to Parshin—Kato theory and integer polytopes,} Funct. Anal. Appl., {\bf 40} (2006), 126--133.

\bibitem{Osipov2003}
D. Osipov, {\em To the multidimensional tame symbol,} Preprints aus dem Institut f\"ur Mathematik, Humboldt-Universit\"at zu Berlin, {\bf 13} (2003), 27 pp., https://doi.org/10.18452/2609; arXiv:1105.1362v1.

\bibitem{OsipovZhu1}
D. Osipov, X. Zhu, {\em A categorical proof of the Parshin reciprocity laws on algebraic surfaces,} Algebra \& Number Theory, {\bf 5} (2011), 289--337.

\bibitem{OsipovZhu2}
D. Osipov, X. Zhu, {\em The two-dimensional {C}ontou-{C}arr{\`e}re symbol and reciprocity laws,} J. Algebr. Geom., {\bf 25} (2016), 703--774.

\bibitem{Par84}
A.~N.~Parshin, {\em Local class field theory,} Trudy Mat. Inst. Steklov., {\bf 165} (1984), 143--170.

\bibitem{Par90}
A.~N.~Parshin, {\em Galois cohomology and the Brauer group of local fields,} Trudy Mat. Inst. Steklov., {\bf 183} (1990), 159--169.

\bibitem{Tate}
J. Tate, {\em Residues of differentials on curves,} Annales scientifiques de l'\'Ecole Normale Sup\'erieure, S{\'e}r.~4, {\bf 1} (1968), 149--159.
	
	
\end{thebibliography}
\end{document}